\documentclass[preprint,12pt,authoryear]{cfg/elsarticle}
\usepackage[letterpaper]{geometry}
\newgeometry{textwidth=6.5in}

\usepackage{amssymb}

\usepackage{amsmath}
\usepackage{graphicx}
\usepackage{enumerate}
\usepackage{url}
\usepackage{array}
\usepackage{subcaption}
\usepackage{amsthm}
\usepackage{fancyhdr}

\usepackage{xcolor}

\newtheorem{definition}{Definition}
\newtheorem{theorem}{Theorem}
\newtheorem{corollary}{Corollary}
\newtheorem{remark}{Remark}

\DeclareMathOperator{\I}{I}

\DeclareMathOperator{\matvec}{vec}
\DeclareMathOperator{\eff}{eff}
\DeclareMathOperator{\E}{E}

\DeclareMathOperator{\D}{D}
\DeclareMathOperator{\Median}{Median}
\DeclareMathOperator{\T}{T}

\journal{arXiv $\qquad\quad\quad\quad$ Copyright \copyright $\,$ 2021 Justin Fishbone}

\begin{document}

\begin{frontmatter}



\title{New Highly Efficient High-Breakdown Estimator of Multivariate Scatter and Location for Elliptical Distributions}

\author[inst1]{Justin A. Fishbone\corref{mycorrespondingauthor}}
\affiliation[inst1]{organization={Bradley Department of Electrical and Computer Engineering\\
    Virginia Polytechnic Institute and State University},
            addressline={\\7054 Haycock Rd}, 
            city={Falls Church},
            postcode={22043}, 
            state={VA},
            country={USA}}
\author[inst1]{Lamine Mili}

\begin{abstract}
High-breakdown-point estimators of multivariate location and shape matrices, such as the MM-estimator with smooth hard rejection and the Rocke S-estimator, are generally designed to have high efficiency at the Gaussian distribution. However, many phenomena are non-Gaussian, and these estimators can therefore have poor efficiency. This paper proposes a new tunable S-estimator, termed the S-q estimator, for the general class of symmetric elliptical distributions, a class containing many common families such as the multivariate Gaussian, t-, Cauchy, Laplace, hyperbolic, and normal inverse Gaussian distributions. Across this class, the S-q estimator is shown to generally provide higher maximum efficiency than other leading high-breakdown estimators while maintaining the maximum breakdown point. Furthermore, its robustness is demonstrated to be on par with these leading estimators while also being more stable with respect to initial conditions. From a practical viewpoint, these properties make the S-q broadly applicable for practitioners. This is demonstrated with an example application---the minimum-variance optimal allocation of financial portfolio investments.
\end{abstract}


\end{frontmatter}


\pagestyle{fancy}
\fancyhf{}
\lhead{Copyright \copyright $\,$ 2021 Justin Fishbone}
\rhead{\thepage}

\pagebreak
\section{Introduction}
\label{sec:intro}

\cite{Huber1964RobustParameter} introduced what is the most common class of estimators, M-estimators. Although originally applied to the location case, \cite{Maronna1976RobustScatter} expanded the definition to include multivariate location and scatter. After the sample median, perhaps the most common robust M-estimators are those using the general rho functions such as the Huber or Tukey bisquare functions.

However, the drawback of using general rho functions is that they have limited efficiency when applied to parameter estimation of nonideal probability distributions. To address this, various M-estimator approaches have been taken to iteratively reweight maximum likelihood estimator (MLE) weights based on the estimated probability density function (PDF) \citep[e.g.,][]{Windham1995RobustifyingFitting,Basu1998RobustDivergence,Choi2000RenderingModel,Ferrari2010MaximumEstimation}.

Even with these improvements, multivariate M-estimators inherently have limited robustness. For example, \cite{Maronna1976RobustScatter} showed that the upper-bound on the breakdown point for $p$-dimensional M-estimators is $(p+1)^{-1}$, which converges to zero with large $p$. To combat this weakness, \cite{Rousseeuw1984RobustS-Estimators} introduced regression S-estimators, which \cite{Davies1987AsymptoticMatrices} expanded to multivariate location and scatter. Davies also showed that the asymptotic breakdown point of S-estimators can be set to 1/2, which is the theoretical maximum of any equivariant estimator.

In practical scenarios however, estimators may have large bias at considerably lower contamination levels than the breakdown point. For many years, the Tukey bisquare was the standard rho function for S-estimators \citep[for example, see][]{Lopuhaa1989OnCovariance,Rocke1996RobustnessDimension}. However, in the context of multivariate S-estimators, the bisquare is not tunable, so its robustness falls off with increasing $p$. For this reason, \cite{Rocke1996RobustnessDimension} introduced the tunable biflat and translated biweight rho functions. \citet[sec.~6.4.4]{Maronna2006RobustMethods} slightly modified the biflat, proposing the Rocke rho function. The Rocke S-estimator (shortened here to S-Rocke) is currently the recommended high-breakdown estimator for large dimensions $(p\geq15)$ \citep[sec.~6.10]{Maronna2017RobustLocation,Maronna2019RobustR}. The recommended estimator for lower dimensions is the MM-estimator with the smoothed hard rejection function (MM-SHR).

There are two major shortcomings of the S-Rocke estimator that will be discussed in this paper. Firstly, it has low efficiency for small dimension, $p$. Although this is an inherent disadvantage of all S-estimators, it is exceptionally acute for the S-Rocke. Secondly, the S-Rocke has poor efficiency for most common non-Gaussian distributions. This is a common problem for general-purpose estimators such as the Rocke and bisquare S-estimators, the MM-SHR, and the Huber and bisquare M-estimators. Examples of common phenomena that are frequently modeled by non-Gaussian distributions include stock returns, radar sea clutter, and speech signals, which approximately follow generalized hyperbolic \citep{KonlackSocgnia2014AReturns}, K- \citep{Ward1990MaritimeSurface}, and Laplace distributions \citep{Gazor2003SpeechDistribution}, respectively.

This paper proposes and explores a new subclass of tunable, maximum-breakdown-point S-estimators that is applicable across common continuous elliptical distributions. This estimator, named the S-q estimator, uses a density-based reweighting to attain generally higher maximum efficiency across the elliptical class as compared to the S-Rocke and MM-SHR estimators. These estimators are compared from the viewpoints of statistical and computational efficiency, robustness, and stability. 

Although the focus on elliptical distributions sounds limiting, as discussed in the next section, most common continuous multivariate distributions---such as the Gaussian, t-, Laplace, and hyperbolic distributions---fall into this class. As \cite{Frahm2009AsymptoticScales} discussed, this assumption is ``fundamental in multivariate analysis.''

This paper is organized as follows. Section \ref{sec:def} defines the new estimator and provides its functions for the most common elliptical distributions. Basic properties related to the consistency of the S-q estimator are summarized in Section \ref{sec:consistency}. Section \ref{sec:efficiency} provides the asymptotic distribution of the S-q estimator and compares the maximum achievable efficiencies of the S-q, S-Rocke, and MM-SHR estimators. In Section \ref{sec:robustness}, the finite-sample breakdown point of the S-q is discussed, the theoretical influence functions of the estimators are compared, and the empirical finite-sample robustness of the estimators are briefly explored. Section \ref{sec:compAnalysis} assesses two computational aspects of the estimators: computational efficiency, and stability with respect to initial estimates. A real-world example in Section \ref{sec:example} demonstrates the application of the estimators for the minimum-variance optimal allocation of financial portfolio investments. Finally, conclusions are summarized in Section \ref{sec:conclusion}.

\section{Defining the S-q Estimator}
\label{sec:def}
This section builds the definition of the proposed S-q estimator. First, the elliptical class of distributions is reviewed. The multivariate S-estimator definition is then summarized, and finally, the S-q is defined.

\subsection{Elliptical Distributions}
\label{subsec:elliptical}
The elliptical distribution is a general class of multivariate probability distributions encompassing many familiar subclasses such as the symmetric Gaussian, t-, Cauchy, Laplace, hyperbolic, variance gamma, and normal inverse Gaussian distributions. Table \ref{tab:distDef} summarizes the most common elliptical distributions (\citealp[][p.~69]{Fang1990SymmetricDistributions}; \citealp[][]{Deng2018OnInequality}).

\begin{table} \caption{Summary of Common Elliptical Distributions} \label{tab:distDef}
\begin{center}
\small
\renewcommand{\arraystretch}{0.8}
\begin{tabular}{l l}
\hline
Distribution Name & Generating Function, $\phi(d)$\\
& \quad\textit{\{Range of Parameters\}}\\
\hline
Kotz type & $d^{N}exp(-rd^s)$ \\
& \quad$\left\{r>0,s>0,\,N>-\frac{p}{2}\right\}$ \\
&\\
Gaussian & $exp\left(-\frac{d}{2}\right)$\\
\textit{(Kotz type with $N=0,s=1,r=1/2$)}&\\
&\\
Pearson type II & $(1-d)^m$, $\quad d\in[0,1]$\\
& \quad$\{m>0\}$ \\
&\\
Pearson type VII & $(1+d/s)^{-N}$\\
& \quad$\left\{N>p/2,\,s>0\right\}$\\
&\\
t & $(1+d/\nu)^{-(\nu+p)/2}$\\
\textit{(Pearson VII with $s=\nu,N=(\nu+p)/2$)} & \quad$\{\nu>0\}$\\
&\\
Cauchy & $(1+d)^{-(1+p)/2}$\\
\textit{(t with $\nu=1$)} & \\
&\\
&\\
Generalized hyperbolic & $\left(\sqrt{\psi(\chi+d)}\right)^{\lambda-p/2}K_{\lambda-p/2}\left(\sqrt{\psi(\chi+d)}\right)$\\
& \quad$\left\{\psi>0,[\chi>0,\lambda\in\mathbb{R} \;or\; \chi=0,\lambda>0]\right\}$\\
&\\
Variance gamma & $\left(\sqrt{\psi\,d}\right)^{\lambda-p/2}K_{\lambda-p/2}\left(\sqrt{\psi\,d}\right)$\\
\textit{(Gen. hyperbolic with $\chi=0$)} & \quad$\left\{\psi>0,\lambda>0\right\}$\\
&\\
Laplace & $\left(\sqrt{2d}\right)^{1-p/2}K_{1-p/2}\left(\sqrt{2d}\right)$\\
\textit{(Variance gamma with $\psi=2,\lambda=1$)} & \\
&\\
Multivariate hyperbolic & $exp\left(-\sqrt{\psi(\chi+d)}\right)$\\
\textit{(Gen. hyperbolic with $\lambda=(p+1)/2$)} & \quad$\left\{\psi>0,\,\chi\geq0\right\}$ \\
&\\
Hyperbolic with univariate marginals & $\left(\sqrt{\psi(\chi+d)}\right)^{1-p/2}K_{1-p/2}\left(\sqrt{\psi(\chi+d)}\right)$ \\
\textit{(Gen. hyperbolic with $\lambda=1$)} & \quad$\left\{\psi>0,\,\chi\geq0\right\}$ \\
&\\
Normal inverse Gaussian & $\left(\sqrt{\psi(\chi+d)}\right)^{-(1+p)/2}K_{-(1+p)/2}\left(\sqrt{\psi(\chi+d)}\right)$ \\
\textit{(Gen. hyperbolic with $\lambda=-1/2,\,\chi>0$)} & \quad$\left\{\psi>0,\,\chi>0\right\}$\\
\hline
\end{tabular}
\end{center}
\end{table}

Symmetric elliptical distributions are defined as being a function of the
squared Mahalanobis distance,\footnote{Some texts define the Mahalanobis distance with the mean and covariance, but this more restrictive definition excludes thick-tailed distributions where these do not exist, such as Cauchy distributions.}
$d\left(\boldsymbol{x},\boldsymbol{\mu},\boldsymbol{\Sigma} \right)=\left(\boldsymbol{x}-\boldsymbol{\mu}\right)^{\T}\boldsymbol{\Sigma}^{-1}\left(\boldsymbol{x}-\boldsymbol{\mu}\right)$,
where $\boldsymbol{x} \in \mathbb{R}^p$, the location 
$\boldsymbol{\mu} \in \mathbb{R}^p$, and the $p\times p$ positive definite symmetric ($\operatorname{PDS}(p)$) scatter 
$\boldsymbol{\Sigma} \in \operatorname{PDS}(p)$. When the PDF is defined, it has the form
$f_X\left(\boldsymbol{x}\right)=\alpha_p|\boldsymbol{\Sigma}|^{-1/2} \phi\left(d\left(\boldsymbol{x},\boldsymbol{\mu},\boldsymbol{\Sigma} \right)\right),$
for some generating function $\phi(d)$, and where $\alpha_p$ is a constant that ensures $f_X\left(\boldsymbol{x}\right)$ integrates to one. Table \ref{tab:distDef} lists common generating functions. When the covariance exists, it is proportional to the scatter matrix, $\boldsymbol{\Sigma}$. The corresponding shape matrix is commonly defined as
\begin{equation} \label{eq:shapeDef}
    \boldsymbol{\Omega}=\boldsymbol{\Sigma}/|\boldsymbol{\Sigma}|^{1/p}.
\end{equation} 
The PDF of $d\left(\boldsymbol{x},\boldsymbol{\mu},\boldsymbol{\Sigma} \right)$ is given by \citep{Kelker1970DistributionGeneralization} 
\begin{equation} \label{eq:f(d)} f_D\left(d\right)=\beta_p\,d^{p/2-1}\phi\left(d\right), \end{equation}
where $\beta_p=\alpha_p\pi^{p/2}/\Gamma(p/2)$. Hereafter, all densities, $f(d),$ refer to the density of $d\left(\boldsymbol{x},\boldsymbol{\mu},\boldsymbol{\Sigma} \right)$ in (\ref{eq:f(d)}), so the subscript $D$ will be omitted. It is also generally assumed that $p>2.$

\subsection{S-Estimators}
\label{subsec:sest}

Given a set of $n$ $p$-dimensional samples, 
$\{\boldsymbol{x}_1, ..., \boldsymbol{x}_n\}$, S-estimators of location and shape are defined as \citep[Sec. 6.4.2]{Maronna2006RobustMethods}
\begin{equation} \label{eq:sEstDef}
\begin{aligned}
\left( \widehat{\boldsymbol{\mu}}, \widehat{\boldsymbol{\Omega}} \right) = \text{arg\,min}
&&& \widehat{\sigma}\\
\text{subject to}
&&& |\boldsymbol{\Omega}|=1, \\
&&& \frac{1}{n} \sum_{i=1}^{n}\rho\left(\frac
        {d\left(\boldsymbol{x}_i,\boldsymbol{\mu},\boldsymbol{\Omega} \right)}
        {\widehat{\sigma}}
    \right)=b,
\end{aligned}
\end{equation}
for some scalar rho function, $\rho(t)$. A \emph{proper} S-estimator rho function should be a continuously differentiable, nondecreasing function in $t\geq0$ with $\rho(0)=0$, and where there is a point $c$ such that $\rho(t)=\rho(\infty)$ for $t \geq c$. For simplicity, and without loss of generality, the rho functions will be normalized so $\rho(\infty)=1$. The parameter $b$ is a scalar that affects the efficiency (see Section \ref{sec:efficiency}) and robustness of the estimator. The purpose of S-estimators is to achieve high robustness, so they are usually configured with $b=1/2-(p+1)/(2n),$ which achieves the maximum theoretical breakdown point that any affine equivariant estimator may have (see Section \ref{subsec:breakdownpoint}). To understand the derivation of the proposed estimator in the next section, note that $\widehat{\sigma}$ in the constraint is an M-estimator of the scale of $d\left(\boldsymbol{\mu},\boldsymbol{\Omega} \right).$ Local solutions of (\ref{eq:sEstDef}) can be found iteratively using the weighted sums $\sum_{i=1}^{n}w\left(d_i/\widehat{\sigma}\right)\left(\boldsymbol{x}_i-\widehat{\boldsymbol{\mu}}\right)=\boldsymbol{0}$ and $\sum_{i=1}^{n}w\left(d_i/\widehat{\sigma}\right)\left(\boldsymbol{x}_i-\widehat{\boldsymbol{\mu}}\right)\left(\boldsymbol{x}_i-\widehat{\boldsymbol{\mu}}\right)^{\T}\propto\widehat{\boldsymbol{\Omega}},$
where the weight function $w(t)=\rho'(t),$ and where $\widehat{\boldsymbol{\Omega}}$ is re-normalized with each iteration. For the empirical results in this paper, the estimators will all be solved using this weighted-sum algorithm.

To estimate the scatter metrix, a separate estimator of $|\boldsymbol{\Sigma}|^{1/p}$ can then be used to scale $\widehat{\boldsymbol{\Omega}}$ using (\ref{eq:shapeDef}). \citet[p.~186]{Maronna2006RobustMethods} discussed a simple estimator to scale $\widehat{\boldsymbol{\Omega}}$ to $\widehat{\boldsymbol{\Sigma}}$. When $\boldsymbol{x}$ is normally distributed, $d$ has a chi-squared distribution with $p$ degrees of freedom. Therefore, they suggested using $\widehat{\boldsymbol{\Sigma}}=\Median\left\{d\left(x_1,\widehat{\boldsymbol{\mu}},\widehat{\boldsymbol{\Omega}}\right), \dots, d\left(x_n,\widehat{\boldsymbol{\mu}},\widehat{\boldsymbol{\Omega}}\right)\right\}\left(\chi^2_p(0.5)\right)^{-1}\widehat{\boldsymbol{\Omega}}$, where $\chi^2_p(0.5)$ is the 50\textsuperscript{th} percentile of the chi-squared distribution. For the general case of elliptical distributions, we propose extending this to 
\begin{equation*}
\widehat{\boldsymbol{\Sigma}}=\frac{\Median\left\{d\left(x_1,\widehat{\boldsymbol{\mu}},\widehat{\boldsymbol{\Omega}}\right), \dots, d\left(x_n,\widehat{\boldsymbol{\mu}},\widehat{\boldsymbol{\Omega}}\right)\right\}}{F^{-1}(0.5)}\widehat{\boldsymbol{\Omega}},
\end{equation*}
where $F(d)$ is the distribution function corresponding to (\ref{eq:f(d)}), and therefore $F^{-1}(0.5)$ is the 50\textsuperscript{th} percentile of the distribution.

For the location and shape matrices, the S-estimator formulation in (\ref{eq:sEstDef}) is equivalent to the alternative one given by 
\begin{equation} \label{eq:sEstALTDef}
\begin{aligned}
\left( \widehat{\boldsymbol{\mu}}, \widehat{\boldsymbol{\Sigma}} \right) = \text{arg\,min}
&&& |\boldsymbol{\Sigma}|\\
\text{subject to}
&&& \frac{1}{n} \sum_{i=1}^{n} \rho\left(\frac{d\left(\boldsymbol{x}_i,\boldsymbol{\mu},\boldsymbol{\Sigma} \right)}{\sigma}\right)=b,
\end{aligned}
\end{equation}
which requires that $\sigma$ be defined such that $b=\E\left[\rho\left(d\left(\boldsymbol{x};\boldsymbol{\mu},\boldsymbol{\Sigma} \right)/\sigma\right)\right]$ for a consistent estimator of $\boldsymbol{\Sigma}$ at an assumed elliptical distribution \citep{Rocke1996RobustnessDimension}. While the first formulation is better for understanding the derivation of the proposed S-q estimator, this second formulation is better for defining and understanding its properties \citep[for example, see][]{Lopuhaa1989OnCovariance}. The scale parameters in the two formulations are related asymptotically by $\sigma=|\boldsymbol{\Sigma}|^{-1/p}\E\left[\widehat{\sigma}\right],$ at the assumed distribution.

The two most common multivariate S-estimators are the bisquare and Rocke \citep[sec.~6.4.2,~6.4.4]{Maronna2019RobustR}. The S-bisquare is given by $\rho_\text{bisq}(t)=\min\left\{1,1-(1-t)^3\right\}$ and $w_\text{bisq}(t)=3(1-t)^2\I(t\leq1),$
which does not have a tuning parameter to control efficiency and robustness. The S-Rocke is given by 
\begin{equation*}
\rho_\gamma(t)=
    \begin{cases}
    0 & \text{if}\;0\leq t \leq 1-\gamma \\
    \frac{t-1}{4\gamma}\left[3-\left(\frac{t-1}{\gamma}\right)^2\right]+\frac{1}{2}
    & \text{if}\;1-\gamma<t<1+\gamma \\
    1 & \text{if}\;1+\gamma \leq t \\
    \end{cases},
\end{equation*}
\[w_\gamma(t)=\frac{3}{4\gamma}\left[1-\left(\frac{t-1}{\gamma}\right)^2\right]\I(1-\gamma \leq t \leq 1+\gamma),\]
where the parameter $\gamma \in (0,1]$ tunes the estimator's efficiency and robustness. The Rocke's maximum efficiency is generally limited at $\gamma=1$, which is extremely restricting for small $p$. Both $\rho_\text{bisq}(t)$ and $\rho_\gamma(t)$ are generic functions that do not depend on the underlying distribution. In the following section, an alternative S-estimator is defined that accounts for the underlying distribution and that generally has better performance across the most common elliptical distributions. It also does not have the same inherent restrictions for small $p$ as $\rho_\gamma(t).$

\subsection{Elliptical Density-Based S-q Estimator}
\label{subsec:def}
The rho function corresponding to the maximum likelihood estimator, $\widehat{\sigma},$ of the scale of $d\left(\boldsymbol{\mu},\boldsymbol{\Omega} \right),$ or equivalently $d\left(\boldsymbol{\mu},\boldsymbol{\Sigma} \right),$ is
$\rho_\text{mle}(t)=-t\,f'\left(t\right)/f\left(t\right).$
We propose weighting this by the power transform of the density,
$\Tilde{\rho}_q\left(t\right)=f(t)^{1-q}\rho_\text{mle}(t),$
where the scalar $q\leq1$ controls the estimator robustness, with $q=1$ corresponding to the maximum likelihood function, and with decreasing $q$ increasing the estimator robustness. In most cases, this rho function is not monotone, as required by S-estimators, so it is denoted with a tilde. This rho function is equivalent to the M-Lq and other M-estimators proposed, for example, by \cite{Windham1995RobustifyingFitting,Basu1998RobustDivergence,Choi2000RenderingModel}; and \cite{Ferrari2010MaximumEstimation}. However, in this particular case of estimating the scale of the squared Mahalanobis distance of an elliptically distributed random vector, the density and rho function do not need to be regenerated with each numerical iteration, $i,$ based on the estimates $\widehat{\boldsymbol{\mu}}^{(i)}$ and $\widehat{\boldsymbol{\Omega}}^{(i)}.$  
Substituting the PDF from (\ref{eq:f(d)}),
\begin{equation} \label{eq:rhoTild}
\Tilde{\rho}_q\left(t\right)=-\left(\beta_p\phi(t)\right)^{s_q}t^{s_p s_q}\left(t\frac{\phi'(t)}{\phi(t)}+s_p\right),
\end{equation}
where $s_p=p/2-1$ and $s_q=1-q$.
Taking the derivative of $\Tilde{\rho}_q(t)$, the corresponding weight function is given by
\begin{equation} \label{eq:wTild}
\Tilde{w}_q(t)=-\left(\beta_p\phi(t)\right)^{s_q}t^{s_p s_q}\left(
\frac{s_q s_p^2}{t} + \left(2 s_q s_p+1\right)\frac{\phi'(t)}{\phi(t)}-qt\left(\frac{\phi'(t)}{\phi(t)}\right)^2+t\frac{\phi''(t)}{\phi(t)}
\right).
\end{equation}
For simplicity, when $q<1$, the scalar $\beta_p$ can be dropped from the calculation of $\Tilde{\rho}_q(t)$ and $\Tilde{w}_q(t)$ in (\ref{eq:rhoTild}) and (\ref{eq:wTild}). When $\phi(t)$ is only positive over a finite domain (e.g. Pearson Type II distribution), then we define $\tilde{\rho}_q(t)$ and $\tilde{w}_q(t)$ to be zero outside this domain.

For the common elliptical distributions listed in Table \ref{tab:distDef}, $\Tilde{\rho_q}(t)$ is monotone in its central region between its global extrema when using appropriate values for $q$ (defined below). The first extremum is the minimum, which we label point $a$, and the second is the maximum, labeled $c$. The distance between $a$ and $c$ varies monotonically with respect to $q$. We use this to define a tunable, double-hard-rejection S-estimator rho function. The value of $\Tilde{\rho}_q(t)$ is held constant between zero and $a$ at value $\Tilde{\rho}_q(a),$ which hard rejects inliers, and the value of $\Tilde{\rho}_q(t)$ is held constant above $c$ at value $\Tilde{\rho}_q(c),$ which hard rejects outliers. The resulting monotonic function is then scaled and shifted so it ranges from zero to one. This defines the S-q estimator.

\begin{definition}
Assuming $\phi(t)$ is twice continuously differentiable over its region of support and 
$\frac{s_q s_p^2}{t} + \left(2 s_q s_p+1\right)\frac{\phi'(t)}{\phi(t)}-qt\left(\frac{\phi'(t)}{\phi(t)}\right)^2+t\frac{\phi''(t)}{\phi(t)}$
has one or two zeros in $t\in(0,\infty)$ for $q<1$, the S-q estimator is the S-estimator with the rho function given by
\begin{equation} \label{eq:rho}
\rho_q(t)=
\begin{cases}
0                                                   & \text{if  $q<1$ and $t \leq a$} \\
s_1\left(\Tilde{\rho}_q(t)-\Tilde{\rho}_q(a)\right) & \text{if $q<1$ and $a<t<c$} \\
1                                                   & \text{if $q<1$ and $t \geq c$} \\
\Tilde{\rho}_q(t)                                   & \text{if $q=1$}
\end{cases},
\end{equation}
where $s_1=\left(\Tilde{\rho}_q(b)-\Tilde{\rho}_q(a)\right)^{-1}$.
The S-q estimator of Type I is the case with one zero (i.e. $a=0$), and the Type II S-q estimator is the case with two zeros.
\end{definition}

For most distributions, $\lim_{q\rightarrow1}c=\infty,$ or at $q=1$, $\Tilde{\rho}_q(t)$ is not bounded. Therefore, we do not scale or shift $\Tilde{\rho}_q(t)$ in this case, and $\rho_q(t)$ is not a proper S-estimator rho function. However, when $q=1$ and $b=1$, the MLE of the scale of $d$ is obtained.
The S-q weight function is the derivative of $\rho_q(t)$ and is given by
\begin{equation} \label{eq:w}
w_q(t)=
\begin{cases}
0                 & \text{if $q<1$ and $t \leq a$} \\
s_1\Tilde{w}_q(t) & \text{if $q<1$ and $a<t<c$} \\
0                 & \text{if $q<1$ and $t \geq c$} \\
\Tilde{w}_q(t)    & \text{if $q=1$}
\end{cases}.
\end{equation}

Table \ref{tab:rhoAC} lists expressions for the \emph{inlier rejection point}, $a,$ and the \emph{outlier rejection point}, $c$,  for the common elliptical distributions in Table \ref{tab:distDef}. For most of these distributions, the equation $\Tilde{w}_q(t)=0$ is quadratic, which provides a closed-form solution for the values of $a$ and $c$. 

\begin{table} \caption{S-q Inlier and Outlier Rejection Points for Common Elliptical Distributions} \label{tab:rhoAC}
\begin{center}
\small
\renewcommand{\arraystretch}{1}
\begin{tabular}{l l p{10cm}}
\hline
Distribution & \multicolumn{2}{l}{Inlier Rejection Point $a$ and Outlier Rejection Point $c$}\\
\hline
Kotz type & $a,c$&$=\left(\frac{s+2s_qN+2s_ps_q\mp\sqrt{s^2+4ss_qN+4ss_ps_q}}{2s_qrs}\right)^{1/s}$ \\
&\\
Gaussian & $a,c$&$=\frac{1+2s_ps_q\mp\sqrt{1+4s_ps_q}}{s_q}$ \\
&\\
Pearson type II & $a,c$&$=\frac{2s_qs_p^2+m\left(2s_qs_p+1\right)\mp\sqrt{m^2\left(4s_qs_p+1\right)+4ms_qs_p^2}}{2\left(s_qs_p^2+m\left(2s_qs_p+ms_q\right)\right)}$ \\
&\\
Pearson type VII & $a,c$&$=s\frac{2Ns_qs_p+N-2s_qs_p^2\mp\sqrt{4N^2s_qs_p-4Ns_qs_p^2+N^2}}{2s_q\left(s_p^2-2Ns_p+N^2\right)}$ \\
&\\
Generalized hyperbolic & $a$&
\(
=
\begin{cases}
0 & when\;\chi=0\;and\;\lambda=1 \\
\left\{t|\Tilde{w}_q(t)=0\;and\;t \in (0,c)\right\} & otherwise \\
\end{cases}
\) \\
& $c$&$=\left\{t|\Tilde{w}_q(t)=0\;and\;t \in (a,\infty)\right\}$\\
\hline
\end{tabular}
\end{center}
\end{table}

The \emph{asymptotic rejection probability} (ARP) is defined as $Pr(d/\widehat{\sigma}\geq c)$ \citep{Rocke1996RobustnessDimension}. Table \ref{tab:rhoAC} can be used to determine $q$ from a desired ARP using $F^{-1}(ARP)$. However, since $w_q(t)$ is very tapered (i.e. applying little weight to values just below $c$), practitioners may choose alternative approaches to tuning that allow for higher estimator efficiencies. For example, the approach used in this paper as well as in \cite{Maronna2017RobustLocation} is to tune the estimators to a desired expected efficiency, which is defined in the next section.

The general definition in (\ref{eq:rho}) specifies that $q\leq1$. In a few particular cases, however, there are some minor restrictions on $q$ (when $q<1$) in order to ensure that $a$ and $c$ are in the support of $f(d)$. Table \ref{tab:rhoQ} lists these restrictions.

\begin{table} \caption{Restrictions on Parameter $q$ for Common Elliptical Distributions} \label{tab:rhoQ}
\begin{center}
\small
\renewcommand{\arraystretch}{1}
\begin{tabular}{l l}
\hline
Distribution & Valid Range of $q$\\
\hline
Kotz type & $q\leq1$ unless $-1-s_p<N<-s_p$, then $1+\frac{s}{4\left(s_p+N\right)}<q\leq1$ \\
Gaussian & $q\leq1$ \\
Pearson type II & $q=1$ or $q<1-\frac{1}{m}$ \\
Pearson type VII & $q\leq1$\\
Generalized hyperbolic\textsuperscript{*} & $q\leq1$ unless $\chi=0$ and $\lambda<1$, then $\textit{unknown}<q\leq1$\\
\hline
\multicolumn{2}{l}{\textit{\textsuperscript{*}Empirically inferred. Computational precision restricts $q\notin(0.998,1),$ approximately.}}\\
\end{tabular}
\end{center}
\end{table}

\begin{figure} 
\begin{center}
\includegraphics[width=6in]{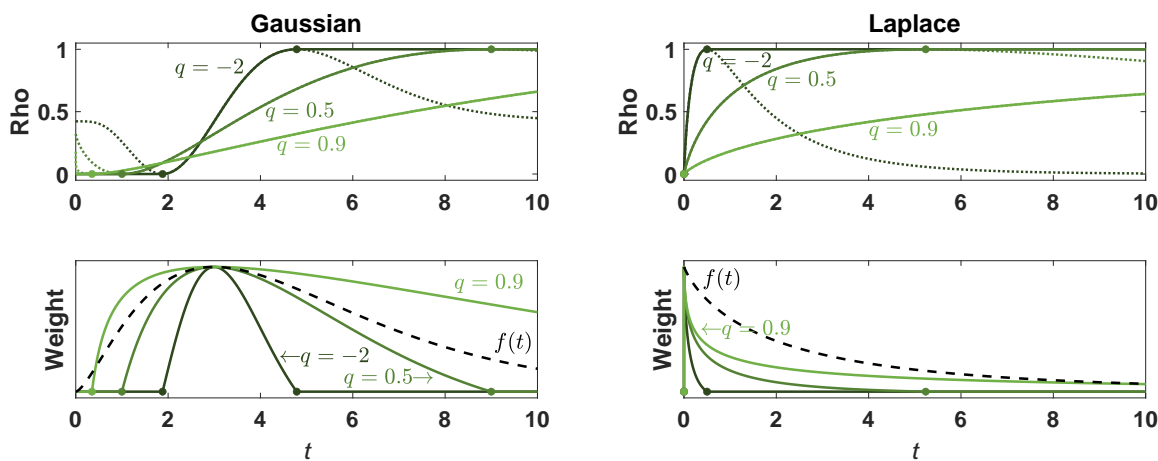}
\end{center}
\caption{Example S-q Rho and Weight Functions for Gaussian (Type II S-q) and Laplace (Type I S-q) Distributions. Rho functions $\rho_q(t)$ (top) and weight functions $w_q(t)$ (bottom) are plotted for the Gaussian distribution (left) and the Laplace distribution (right) for $q \in \{-2, 0.5, 0.9\}$ and $p=5$. On the top, the dotted lines depict the corresponding $\Tilde{\rho}_q(t)$ functions, scaled and shifted to match $\rho_q(t)$. On the bottom, the dashed line depicts the density function. The solid dots indicate points $a,$ when $\Tilde{\rho}_q(t)=0,$ and $c,$ when $\Tilde{\rho}_q(t)=1$. \label{fig:exRhosWs}}
\end{figure}

Figure \ref{fig:exRhosWs} illustrates examples of the S-q functions $\Tilde{\rho}_q(t)$, $\rho_q(t)$, and $w_q(t)$ for the five-dimensional Gaussian (S-q Type II) and Laplace (S-q Type I) distributions and for various values of $q$. As $q$ is decreased, the region of positive weights (area between points $a$ and $c$) narrows, corresponding to increased robustness. The PDF is also plotted, illustrating how $w_q(t)$ roughly follows $f(t)$ in the central region.

\begin{figure} 
\begin{center}
\includegraphics[width=3in]{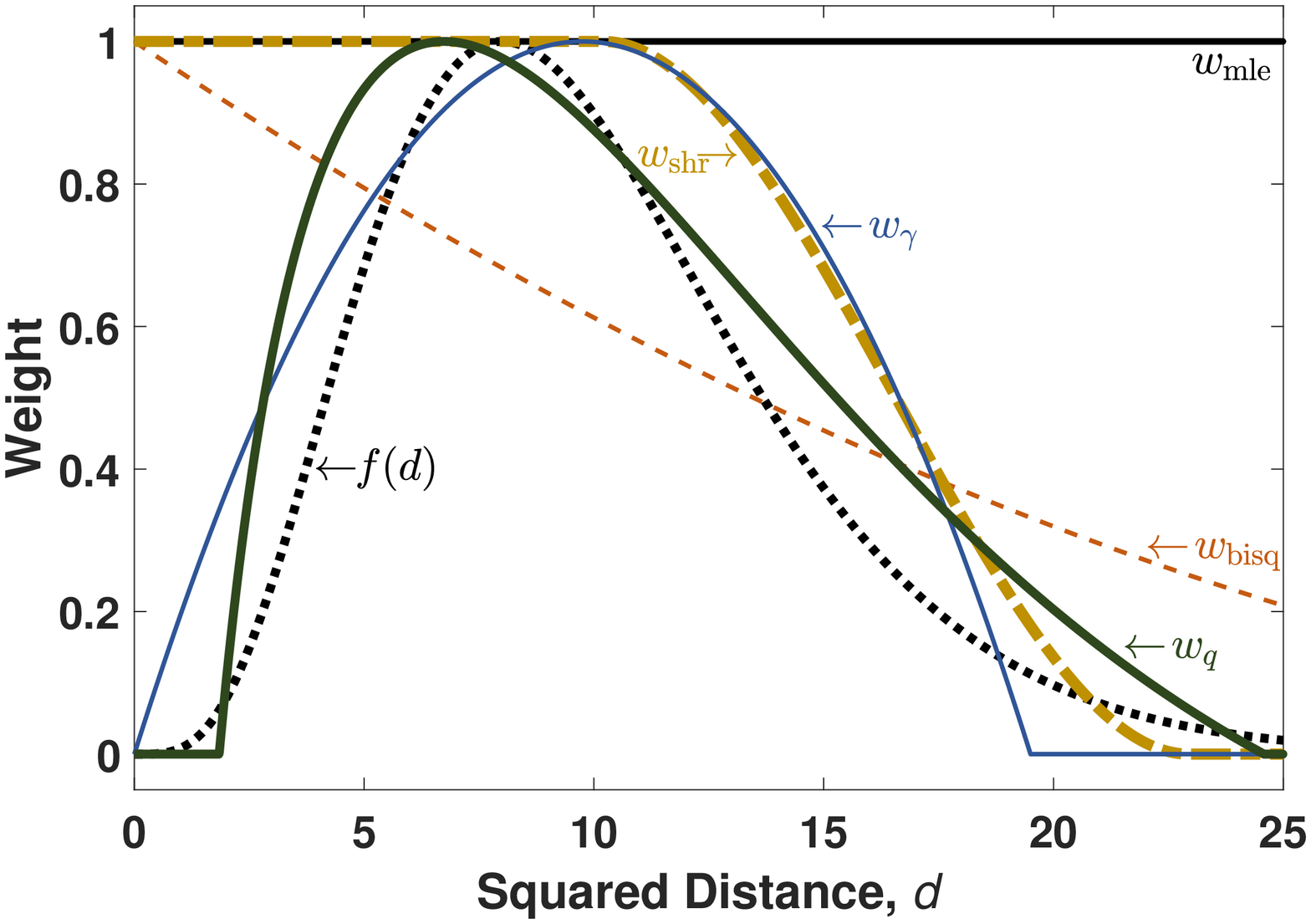}
\end{center}
\caption{Example Comparison of Weight Functions for Various Estimators. For the 10-dimensional Gaussian distribution, the plot depicts the asymptotic weights for the S-q ($w_q$) and MM-SHR ($w_\text{shr}$) estimators tuned to 80\% asymptotic relative efficiency, the S-Rocke ($w_\gamma$) estimator tuned to its maximum efficiency (77\% for this case), and the non-tunable MLE ($w_\text{mle}$) and S-bisquare ($w_\text{bisq}$) estimators. The estimators are set to their maximum breakdown points.
\label{fig:compareWs}}
\end{figure}

Figure \ref{fig:compareWs} compares the S-q asymptotic weights with those of the MM-SHR, S-Rocke, S-bisquare, and maximum likelihood estimators, and with the corresponding PDF. The underlying model is a 10-dimensional standard Gaussian distribution. The MM-SHR and S-q estimators have been tuned to 80\% asymptotic efficiency relative to the MLE. The S-Rocke estimator is tuned to its maximum efficiency, which is 77\% in this instance. The estimators have been set to the maximum breakdown point, with $b=1/2,$ which results in the shifts of the peaks of the weight curves relative to the PDF.

From the figure, it is clear that the Gaussian MLE (i.e. sample estimator) gives uniform weight to all samples, no matter how improbable. The S-Rocke has a quadratic weight function, which is a hard cutoff that cannot capture the tails of $f(d)$. The SHR weight function is cubic, and its shape better captures the shape of the right-half of the PDF. However, the SHR function is designed to approximate a step function, which is poorly suited for many distributions (c.f. $w_{shr}(t)$ in Figure \ref{fig:compareWs} with the Laplace $f(d)$ in Figure \ref{fig:exRhosWs}). Only the S-q weight function follows the general shape of the PDF---giving less weight to less probable observations.

\section{Consistency Properties of the S-q Estimator}
\label{sec:consistency}
As an S-estimator, the S-q estimator inherits properties from its parent class, such as affine equivariance. This section briefly summarizes properties related to its consistency. For more detailed discussion on these, see \citep{Davies1987AsymptoticMatrices}. Here, we use the alternative S-estimator formulation given by (\ref{eq:sEstALTDef}) under the assumptions (A1) that 
\begin{align*}
\text{(A1)\quad\quad} &
b=\E[\rho_q\left(d\left(\boldsymbol{x}_i,\boldsymbol{\mu},\boldsymbol{\Sigma} \right)/\sigma \right)],
\\&\boldsymbol{x}\sim f_X\left(\boldsymbol{x},\boldsymbol{\mu},\boldsymbol{\Sigma},\phi(d)\right)\text{, where $\boldsymbol{x}_i$ are i.i.d.,}
\\&\text{$\phi(d)$ is non-increasing, and}
\\&\text{$\phi(d)$ and $-\rho_q(d/\sigma )$ have common point(s) of decrease.}
\end{align*}

\begin{theorem}[Uniqueness]
Given (A1), minimizing $|\widehat{\boldsymbol{\Sigma}}|$ subject to 
\begin{equation*}
\int_0^\infty \rho_q\left(\frac{d\left(\boldsymbol{x}_i,\widehat{\boldsymbol{\mu}},\widehat{\boldsymbol{\Sigma}} \right)}{\sigma }\right) f_X(d\left(\boldsymbol{x}_i,\boldsymbol{\mu},\boldsymbol{\Sigma} \right))\operatorname{d}\!\boldsymbol{x}=b
\end{equation*}
has a unique solution $\left(\widehat{\boldsymbol{\mu}},\widehat{\boldsymbol{\Sigma}}\right)=\left(\boldsymbol{\mu},\boldsymbol{\Sigma}\right).$
\end{theorem}
\begin{proof}
See \cite[Th. 1]{Davies1987AsymptoticMatrices}.
\end{proof}

\begin{theorem}[Existence]
Given (A1) and $n\geq(p+1)/(1-b/\rho_q(\infty)),$ then the S-q estimator has at least one solution with probability one.
\end{theorem}
\begin{proof}
See \cite[Th. 2]{Davies1987AsymptoticMatrices}.
\end{proof}

\begin{theorem}[Consistency] \label{thm:consistency}
Given (A1), $b=\E[\rho_q(d/\sigma )],$ and $p+1\leq n(1-b/\rho_q(\infty)),$ then \[\lim_{n\rightarrow\infty}\left(\widehat{\boldsymbol{\mu}}_n,\widehat{\boldsymbol{\Sigma}}_n\right)=\left(\boldsymbol{\mu},\boldsymbol{\Sigma}\right).\]
\end{theorem}
\begin{proof}
See \cite[Th. 3]{Davies1987AsymptoticMatrices}. 
\end{proof}

\section{Asymptotic Distribution and Relative Efficiencies}
\label{sec:efficiency}
In this section, the asymptotic distribution of the S-q estimator is provided. From this, measures of efficiency are then defined. Finally, the efficiency of the S-q estimator is compared with leading high-breakdown point estimators.

\subsection{Asymptotic Distribution}
For the asymptotic distribution of the S-q estimate, we continue to use the alternative S-estimator formulation given by (\ref{eq:sEstALTDef}). \cite{Lopuhaa1997AsymptoticCovariance} derived the distribution of S-estimators with assumptions appropriate for the S-q estimator, that is
\begin{align*}
\text{(A2)\;\;\;\;\;\;}&\text{$\phi'_p(t)$ is decreasing with $\phi'_p(d)<0$.}
\end{align*}

Here, we use the following notation. The matrix $\boldsymbol{I}_{p^2}$ is the $p^2\times p^2$ identity matrix, $\boldsymbol{K}_{p^2}$ is the $p^2\times p^2$ commutation matrix, $\otimes$ is the Kronecker product operator, and the operator $\matvec\left(\boldsymbol{\Sigma}\right)$ stacks the columns if $\boldsymbol{\Sigma}$ into a column vector.

\begin{theorem}[Asymptotic distribution] \label{th:asymptoticnormality}
Given (A1) and (A2), the asymptotic distribution of the S-q estimate of
$\left(\widehat{\boldsymbol{\mu}}_n, \widehat{\boldsymbol{\Sigma}}_n\right)$ 
is given by $\sqrt{n}\left(\widehat{\boldsymbol{\mu}}_n-\boldsymbol{\mu}, \widehat{\boldsymbol{\Sigma}}_n-\boldsymbol{\Sigma}\right)\overset{d}{\rightarrow}\left(\boldsymbol{a},\boldsymbol{B}\right)$,  
with $\boldsymbol{a}\perp\boldsymbol{B}$.  
The vector $\boldsymbol{a}\sim\mathcal{N}\left(\boldsymbol{0},\boldsymbol{\Gamma}_{\boldsymbol{\mu}}\right)$  where 
\begin{equation} \label{eq:gammaMu}
    \boldsymbol{\Gamma}_{\boldsymbol{\mu}}=\frac{\omega_1}{\omega_2^{2}}\boldsymbol{\Sigma},
\end{equation}
with 
$\omega_1=p^{-1}\E\left[dw_q^2(d/\sigma )\right]$  
and 
$\omega_2=-2\beta\int_0^\infty p^{-1}d^{p/2}w_q(d/\sigma )\phi'(d)\;\operatorname{d}\!d$. 
The matrix $\boldsymbol{B}\sim\mathcal{N}\left(\boldsymbol{0},\boldsymbol{\Gamma}_{\boldsymbol{\Sigma}}\right)$ 
where 
\begin{equation} \label{eq:gammaSigma}
\boldsymbol{\Gamma}_{\boldsymbol{\Sigma}}=\zeta_1\left(\boldsymbol{I}_{p^2}+\boldsymbol{K}_{p^2}\right)
    \left(\boldsymbol{\Sigma}\otimes\boldsymbol{\Sigma}\right)
    +\zeta_2\matvec\left(\boldsymbol{\Sigma}\right)\matvec\left(\boldsymbol{\Sigma}\right)^{\T},
\end{equation}
with
$\zeta_1=\lambda_1^{-2}p(p+2)\E\left[(d/\sigma )^2w_q^2(d/\sigma ) \right]$   
and
$\zeta_2=\lambda_2^{-2}\E\left[ (\rho_q(d/\sigma )-b)^2\right] - 2p^{-1}\zeta_1$, 
where
$\lambda_1=-2\beta\int_0^\infty \sigma^{-1}d^{p/2+1}w_q(d/\sigma )\phi'(d)\;\operatorname{d}\!d$ 
and 
$\lambda_2=-\beta\int_0^\infty d^{p/2}\left(\rho_q(d/\sigma )-b\right)\phi'(d)\;\operatorname{d}\!d.$
\end{theorem}
\begin{proof}
See \cite[Corollary 2]{Lopuhaa1997AsymptoticCovariance}.
\end{proof}

\cite{Frahm2009AsymptoticScales} derived the asymptotic distribution of shape matrix estimates for affine equivariant estimators. This enables us to state the asymptotic distribution of the S-q shape estimate, which is applicable using either S-estimator formulation, (\ref{eq:sEstDef}) or (\ref{eq:sEstALTDef}).

\begin{theorem}[Shape asymptotic distribution] \label{th:shapeasymptoticnormality}
Given (A1) and (A2), the asymptotic distribution of the S-q estimate of
shape is given by 
$\sqrt{n}\left(\widehat{\boldsymbol{\Omega}}_n-\boldsymbol{\Omega}\right)\sim\mathcal{N}\left(\boldsymbol{0},\boldsymbol{\Gamma}_{\boldsymbol{\Omega}}\right)$ 
where 
\begin{equation}  \label{eq:gammaOmega}
\boldsymbol{\Gamma}_{\boldsymbol{\Omega}}=\zeta_1\left(\boldsymbol{I}_{p^2}+\boldsymbol{K}_{p^2}\right)
    \left(\boldsymbol{\Omega}\otimes\boldsymbol{\Omega}\right)
    -\frac{2\zeta_1}{p}\matvec\left(\boldsymbol{\Omega}\right)\matvec\left(\boldsymbol{\Omega}\right)^{\T},
\end{equation}
with $\zeta_1$ defined as in Theorem \ref{th:asymptoticnormality}.
\end{theorem}
\begin{proof}
See \cite[Corollary 1]{Frahm2009AsymptoticScales}.
\end{proof}

\subsection{Measures of Efficiency}
The asymptotic efficiency of an estimator, at an assumed distribution, is defined as the ratio of the asymptotic variance of the maximum likelihood estimate to the variance of the estimator under consideration. For multivariate estimation, this definition of efficiency is of large dimension---${p\times p}$ for location and ${p^2\times p^2}$ for shape and scatter. However, for affine equivariant estimation of location and scatter of elliptical distributions, the covariance of the estimate depends only on a scalar. Specifically, (\ref{eq:gammaMu}), (\ref{eq:gammaSigma}), and (\ref{eq:gammaOmega}) are general, with only the scalars $\omega_1/\omega_2^2$ \citep{Bilodeau1999TheoryStatistics}, and $\zeta_1$ and $\zeta_2$ \citep{Tyler1982RadialSphericity} depending on the estimator and the generating function $\phi(d)$. Therefore, the asymptotic efficiency of the estimate $\widehat{\boldsymbol{\mu}}$ can be alternatively defined as
\begin{equation*}
    \eff_\infty\left(\widehat{\boldsymbol{\mu}}\right) 
    = \frac{\omega_{1\text{,mle}}/\omega_{2\text{,mle}}^2}
    {\omega_{1,\widehat{\boldsymbol{\mu}}}/\omega_{2,\widehat{\boldsymbol{\mu}}}^2},
\end{equation*}
and the asymptotic efficiency of the estimate $\widehat{\boldsymbol{\Omega}}$ can alternatively be defined as
\begin{equation} \label{eq:asymEff}
    \eff_\infty\left(\widehat{\boldsymbol{\Omega}}\right) = \frac{\zeta_{1,\text{mle}}}{\zeta_{1,\widehat{\boldsymbol{\Omega}}}}.
\end{equation}
It is common to define asymptotic efficiency this way \citep[for example, see][]{Tyler1983RobustnessMatrices,Frahm2009AsymptoticScales}.

Comparing the S-q estimator's efficiency to another estimator can likewise be achieved analytically using, for example, $\zeta_{1,\gamma}/\zeta_{1,q}$ for the S-Rocke estimator, which when the quotient is greater than one, indicates that the S-q has higher asymptotic efficiency than the S-Rocke estimator. For other S-estimators, the asymptotic distribution parameters $\omega_1,$ $\omega_2,$ and $\zeta_1$ are calculated the same as in Theorems \ref{th:asymptoticnormality} and \ref{th:shapeasymptoticnormality} but using their respective weight functions. MM-estimators have the same asymptotic variance and influence function as S-estimators \citep{Rousseeuw2013High-BreakdownScatter}. For MM-estimators, however, $\sigma$ is effectively the tuning parameter, and it can be set accordingly.

In general, finite-sample performance measures are difficult to derive analytically. Instead, it is common to characterize finite-sample performance by empirically characterizing the behavior of metrics derived from the Kullback-Leibler divergence between the estimated and true distribution \citep[for example, see][]{Huang2006CovarianceLikelihood,Ferrari2010MaximumEstimation}. For t-distributions, which includes the Gaussian distribution, the Kullback-Leibler divergence between $t_{\nu}\left(\boldsymbol{\mu},\boldsymbol{\Sigma}\right)$ and $t_{\nu}\left(\widehat{\boldsymbol{\mu}},\widehat{\boldsymbol{\Sigma}}\right)$ is given by \citep{Abusev2015OnAnalysis} 

\[\operatorname{D}\left(\boldsymbol{\mu},\boldsymbol{\Sigma};\widehat{\boldsymbol{\mu}},\widehat{\boldsymbol{\Sigma}}\right)=\frac{1}{2}\left(\operatorname{Tr}\left(\boldsymbol{\Sigma}^{-1}\widehat{\boldsymbol{\Sigma}}\right) + \left(\boldsymbol{\mu}-\widehat{\boldsymbol{\mu}}\right)^T\boldsymbol{\Sigma}^{-1}\left(\boldsymbol{\mu}-\widehat{\boldsymbol{\mu}}\right) - p - \log\left(\frac{|\widehat{\boldsymbol{\Sigma}}|}{|\boldsymbol{\Sigma}|}\right) \right).\]
Following \cite{Maronna2017RobustLocation}, we then define the joint location and scatter finite-sample relative efficiency as
\begin{equation*}
    \eff_n\left(\widehat{\boldsymbol{\mu}},\widehat{\boldsymbol{\Sigma}};\widehat{\boldsymbol{\mu}}_{\text{mle}},\widehat{\boldsymbol{\Sigma}}_{\text{mle}}\right) = 
    \frac{\E\left[\D\left(\boldsymbol{\mu},\boldsymbol{\Sigma};\widehat{\boldsymbol{\mu}}_{\text{mle}},\widehat{\boldsymbol{\Sigma}}_{\text{mle}}\right)\right]}
         {\E\left[\D\left(\boldsymbol{\mu},\boldsymbol{\Sigma};\widehat{\boldsymbol{\mu}},\widehat{\boldsymbol{\Sigma}}\right)\right]},
\end{equation*}
where $\widehat{\boldsymbol{\mu}}_\text{mle}$ and  $\widehat{\boldsymbol{\Sigma}}_\text{mle}$ are the location and scatter matrices corresponding to the maximum likelihood estimate, and where the expectation is calculated empirically using the sample mean over $m$ Monte Carlo trials. The location and the scatter finite-sample relative efficiencies are then respectively defined as 
$\eff_n\left(\widehat{\boldsymbol{\mu}},\boldsymbol{\Sigma};\widehat{\boldsymbol{\mu}}_{\text{mle}},\boldsymbol{\Sigma}\right)$
and
$\eff_n\left(\boldsymbol{\mu},\widehat{\boldsymbol{\Sigma}};\boldsymbol{\mu},\widehat{\boldsymbol{\Sigma}}_{\text{mle}}\right).$
Likewise, we define the shape matrix finite-sample relative efficiency as
\begin{equation} \label{eq:finiteEff}
    \eff_n\left(\widehat{\boldsymbol{\Omega}};\widehat{\boldsymbol{\Omega}}_{\text{mle}}\right) = 
    \frac{\E\left[\D\left(\boldsymbol{\mu},\boldsymbol{\Omega};\boldsymbol{\mu},\widehat{\boldsymbol{\Omega}}_{\text{mle}}\right)\right]}
         {\E\left[\D\left(\boldsymbol{\mu},\boldsymbol{\Omega};\boldsymbol{\mu},\widehat{\boldsymbol{\Omega}}\right)\right]}.
\end{equation}

\subsection{Comparison of Estimator Efficiency}
\label{subsec:effResults}
Any estimator must provide a good estimate in the absence of contamination and when tuned to its maximum efficiency. This section compares the maximum achievable efficiencies of the S-q, S-Rocke, and MM-SHR estimators when set to their maximum breakdown point. The results below cover large swaths of the most common elliptical families in Table \ref{tab:distDef} for a moderate dimension of $p=20.$ These swaths were specifically chosen to cover everyday distributions: Gaussian, Cauchy, Laplace, hyperbolic, and normal inverse Gaussian distributions.

Robust scatter matrix estimation is generally ``more difficult'' than the estimation of location \citep{Maronna2019RobustR}, and as \cite{Maronna2017RobustLocation} demonstrated, divergence and efficiency metrics for scatter matrix estimators are generally much worse than for the corresponding estimators of location. Likewise, due to the high dimensionality of the estimate, the underlying shape matrix is the most difficult part of estimating the scatter matrix. Additionally, many practical applications such as
multivariate regression, principal components analysis, linear discriminant analysis, and canonical correlation analysis
only require the shape matrix, and not the full scatter or covariance matrices \citep{Frahm2009AsymptoticScales}.
Therefore, unless otherwise noted, the performance results in this paper are for the shape matrix, with metrics given by (\ref{eq:asymEff}), (\ref{eq:finiteEff}), and $D\left(\boldsymbol{\mu},\boldsymbol{\Omega};\boldsymbol{\mu},\widehat{\boldsymbol{\Omega}}\right).$

\begin{figure} 
\begin{center}
\includegraphics[width=3in]{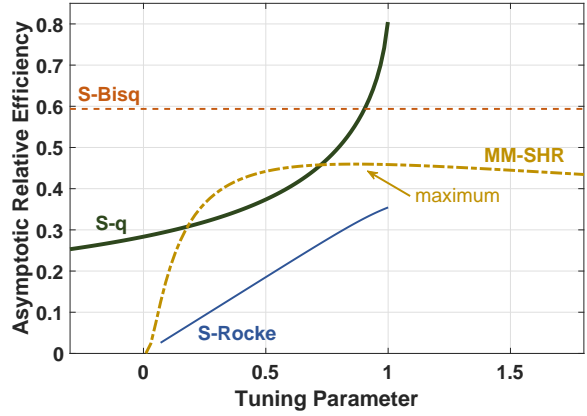}
\end{center}
\caption{Estimator Asymptotic Relative Efficiency versus Tuning Parameter. The relative efficiencies of the estimators are plotted as a function of tuning parameter for the Cauchy distribution with $p=20.$ All estimators are set to the maximum breakdown point. The S-Rocke parameter is in $[0,1]$, the MM-SHR parameter is in $(0,\infty)$, and the S-q parameter is in $(-\infty,1)$ for the maximum breakdown point.
\label{fig:effVsTune}}
\end{figure}

The maximum efficiencies of the S-q and S-Rocke generally occur when their parameters $q$ and $\gamma$ are set to one---although the maximum breakdown point of the S-q is only achieved when $q<1$. However, the maximum efficiency of the MM-SHR must be determined by a search as depicted in Figure \ref{fig:effVsTune}, which plots, as an example, asymptotic efficiency versus tuning parameter for the estimators for the 20-dimensional Cauchy distribution. At the limit, as the MM-SHR parameter is increased toward infinity, all samples receive equal weight, which is the MLE for the Gaussian distribution, but not for distributions such as the Cauchy. In general, for each tunable estimator, its efficiency decreases while its robustness increases as its parameter is decreased. At the lower limit of its parameter, its weight function is a delta function that may reject all the samples and may result in zero efficiency. At this point, the robustness is high, but the weighted-sum solution depends entirely on the initial estimates $\widehat{\boldsymbol{\mu}}^{(0)}$ and $\widehat{\boldsymbol{\Omega}}^{(0)}.$

It should be noted that although generally of high efficiency, the S-q estimate at its limit with $q=1$ is not necessarily the maximum likelihood estimate for location and scatter. The MLE weight function for location and scatter is given by \citep{Tyler1982RadialSphericity}
\begin{equation} \label{eq:mleWeight}
    w_{\text{mle}}(t)=-2\frac{\phi'(t)}{\phi(t)}
\end{equation}
whereas at $q=1,$ (\ref{eq:w}) gives
\begin{equation}  \label{eq:sq1Weight}
w_{q=1}(t)=
-\frac{\phi'(t)}{\phi(t)}+t\left(\frac{\phi'(t)}{\phi(t)}\right)^2-t\frac{\phi''(t)}{\phi(t)}.
\end{equation}

\begin{theorem}[Relation to MLE efficiency]
Assuming $b=\E\left[\rho_q(d\left(\boldsymbol{x};\boldsymbol{\mu},\boldsymbol{\Sigma} \right))\right],$ the asymptotic S-q estimate with $q=1$ is the maximum likelihood estimate for the location and scatter matrices for distributions where
\begin{equation} \label{eq:mleCriteria}
    t\frac{\phi''(t)}{\phi(t)} - t\left(\frac{\phi'(t)}{\phi(t)}\right)^2
    = y\frac{\phi'(t)}{\phi(t)},
\end{equation}
for some value $y.$ Therefore, the S-q estimator can asymptotically achieve the Cramér–Rao lower bound for these distributions.
\end{theorem}
\begin{proof}
The S-estimator scaling of $b=\E\left[\rho_q(d\left(\boldsymbol{x};\boldsymbol{\mu},\boldsymbol{\Sigma} \right))\right]$ results in an estimate that is invariant to scaling of the weight function. Therefore, this theorem follows directly from proportionally equating (\ref{eq:mleWeight}) to (\ref{eq:sq1Weight}).
\end{proof}
\begin{remark}
Although this theorem inherently assumes the alternative S-estimator formulation given by (\ref{eq:sEstALTDef}), it still holds true for location and shape matrices using the primary S-estimator formulation in (\ref{eq:sEstDef}).
\end{remark}
\begin{remark}
If $\lim_{q\rightarrow1}\tilde{\rho}_q(c)$ is finite, and $b\leq1/2,$ then both the Cramér–Rao lower bound (maximum efficiency) and the maximum breakdown point can occur in the limit as $q\rightarrow1$ (see Corollary \ref{th:maxBreakdown} in the next section).
\end{remark}

\begin{figure} 
\begin{center}
\includegraphics[width=3in]{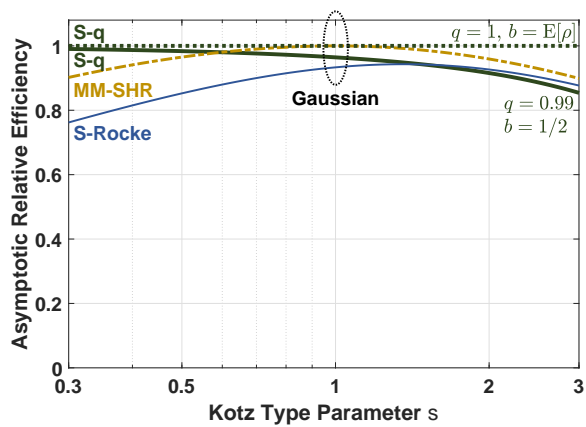}
\end{center}
\caption{Estimator Maximum Achievable Asymptotic Efficiency for Kotz Type Distribution versus Parameter $s.$ Maximum achievable asymptotic shape efficiencies for the maximum breakdown point are plotted for parameters $N=0,$ and $r=1/2.$ The maximum absolute asymptotic shape efficiency of the S-q estimator for $q=1$ is also shown.
\label{fig:maxEff-Kotz}}
\end{figure}

An example family that satisfies this theorem is the Kotz type with parameter $N=0.$ Note, however, that $\lim_{q\rightarrow1}\tilde{\rho}_q(c)=\infty,$ so high breakdown cannot be achieved simultaneously. This is illustrated in Figure \ref{fig:maxEff-Kotz}, which provides the estimators' maximum achievable asymptotic shape efficiencies for the Kotz type distribution with parameters $N=0$ and $r=1/2$ as a function of parameter $s$. In this example, the S-q efficiency is plotted for its maximum absolute efficiency with $q=1$ and for it approximate maximum high-breakdown efficiency with $q=0.99.$ As seen in the figure, the high cost of high-breakdown is particularly acute for large $s.$ 

The remainder of this paper will focus on maximum efficiency at the maximum breakdown point. Figure \ref{fig:maxEff-Kotz} also provides the S-Rocke and MM-SHR estimator's maximum efficiencies at their maximum breakdown points. The MM-SHR efficiency peaks at $s=1,$ which is expected since this is the Gaussian distribution, and the S-Rocke efficiency peaks just above this point. Their efficiencies fall off precipitously for larger and smaller values of $s.$ The efficiency of the S-q, conversely, increases toward unity for smaller $s.$

\begin{figure} 
\begin{center}
\includegraphics[width=6in]{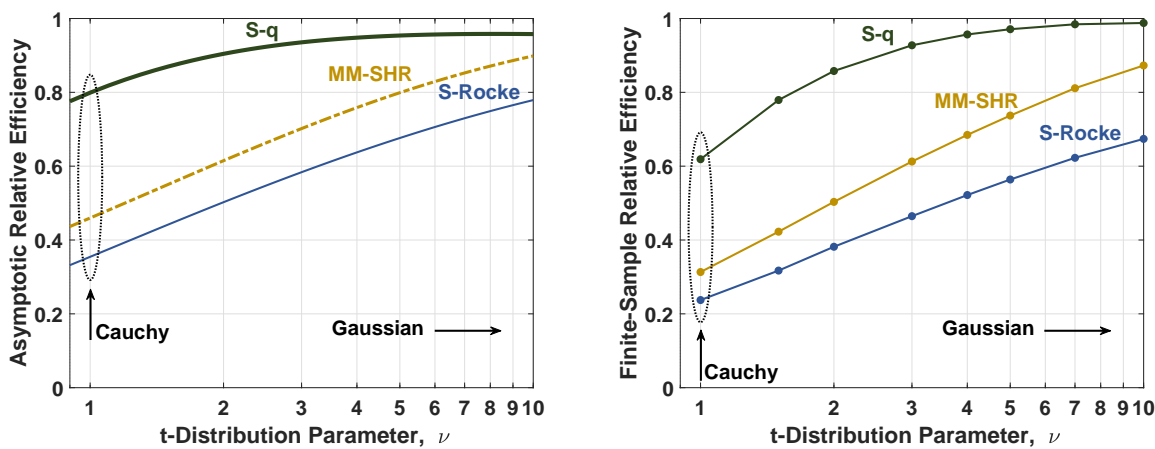}
\end{center}
\caption{Estimator Maximum Achievable Efficiency for t-Distribution versus Distribution Parameter, $\nu.$ Maximum achievable asymptotic (left) and small-sample ($n=3p$; right) efficiencies for the maximum breakdown point are plotted.
\label{fig:maxEff-t}}
\end{figure}

The estimators' maximum achievable asymptotic efficiencies for the t-distribution as a function of the distribution parameter, $\nu$, are plotted on the left of Figure \ref{fig:maxEff-t}. When $\nu=1$, the t-distribution corresponds to a Cauchy distribution, and when $\nu\rightarrow\infty,$ it corresponds to the Gaussian distribution. The S-q estimator offers the highest efficiency of the three estimators for thicker tails.

The maximum achievable small-sample relative efficiencies using $n=3p$ are plotted on the right of Figure \ref{fig:maxEff-t}. The initial estimates were made using the \cite{Pena2007CombiningData} kurtosis plus specific directions (KSD) estimator as recommended and provided by \cite{Maronna2017RobustLocation}. Comparing these finite-sample results with the asymptotic ones on the left, it is seen that the relative results are similar. This general similarity implies that the relative performance of the asymptotic efficiencies can often be a good surrogate for the relative performance of the finite-sample efficiencies when there is no closed-form expression for the divergence in (\ref{eq:finiteEff}).

\begin{figure} 
\begin{center}
\includegraphics[width=6in]{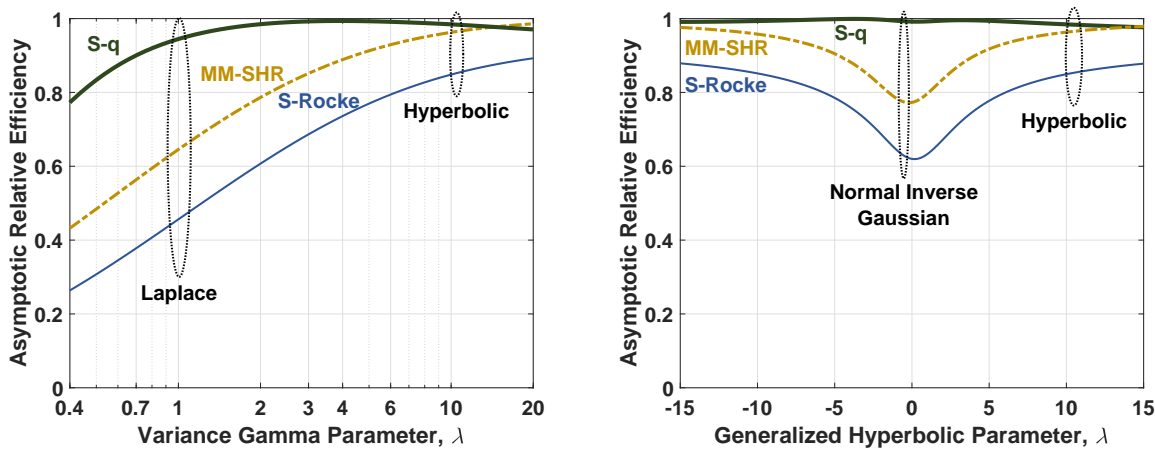}
\end{center}
\caption{Estimator Maximum Achievable Asymptotic Efficiency for Generalized Hyperbolic Distribution versus Parameter $\lambda.$ Maximum achievable asymptotic efficiencies for the maximum breakdown point are plotted for the variance gamma distribution with parameter $\psi=2$ (left) and for the generalized hyperbolic distribution with parameters $\psi=2$ and $\chi=1$ (right).
\label{fig:maxEff-genHyper}}
\end{figure}

The estimators' maximum achievable asymptotic efficiencies for the variance gamma distribution with $\psi=2$ are plotted as a function of parameter $\lambda$ on the left of Figure \ref{fig:maxEff-genHyper}. The plots highlight the Laplace $(\lambda=1)$ and multivariate hyperbolic $(\lambda=(p+1)/2)$ distributions. The S-q exhibits good performance for the hyperbolic and remarkably good performance for the Laplace.

The estimators' maximum achievable asymptotic efficiencies for the generalized hyperbolic distribution with $\psi=2$ and $\chi=1$ are plotted as a function of parameter $\lambda$ on the right of Figure \ref{fig:maxEff-genHyper}. The plots highlight the normal inverse Gaussian $(\lambda=-1/2)$ and hyperbolic $(\lambda=(p+1)/2)$ distributions. The S-q again exhibits good performance for the hyperbolic, and it exhibits remarkably good performance for the normal inverse Gaussian.

\section{Robustness Analysis}
\label{sec:robustness}
The robustness of the S-q estimator is now explored. First, the breakdown point is provided. The influence function is then explored. Finally, finite-sample simulation results are provided to further illustrate the robustness of the high-breakdown estimators.

\subsection{Breakdown Point} \label{subsec:breakdownpoint}
The finite-sample breakdown point of a multivariate estimator of location or scatter is defined as the fraction of the samples, $\epsilon n,$ that can be set to either drive $\|\widehat{\boldsymbol{\mu}}\|=\infty$ or drive an eigenvalue of $\widehat{\boldsymbol{\Sigma}}$ to either zero or infinity. Unlike multivariate M-estimators, which only achieve an asymptotic breakdown point of $(p+1)^{-1}$ \citep{Maronna1976RobustScatter}, S-estimators are able to achieve the maximum possible finite-sample breakdown point that any affine equivariant estimator may have \cite[Th. 6]{Davies1987AsymptoticMatrices}. For the following theorem, the term samples in \emph{general position} means that no more than $p$ samples are contained in any hyperplane of dimension less than $p$.

\begin{theorem}[Finite-sample breakdown point] \label{th:breakdown}
Assuming (A1) and $q<1,$ when $n$ samples are in general position and $n(1-2b)\geq p+1,$ the breakdown point of the S-q estimator is $(\lfloor nb\rfloor+1)/n.$
\end{theorem}
\begin{proof}
As discussed above Section \ref{subsec:def}, $q<1$ ensures a proper S-estimator with finite value $c$ and bounded rho function. See \cite[Th. 5]{Davies1987AsymptoticMatrices}. 
\end{proof}
\begin{corollary} \label{th:maxBreakdown}
The maximum breakdown point is $\lfloor (n-p+1)/2\rfloor/n,$ which is achieved when $b=1/2-(p+1)/(2n).$ Asymptotically, this is $1/2$ at $b=1/2.$
\end{corollary}

\subsection{Influence Function}
The influence function (IF) of an estimator characterizes its sensitivity to an infinitesimal point contamination at $\boldsymbol{z}\in\mathbb{R}^p,$ standardized by the mass of the contamination, $\epsilon.$ The influence function for estimator $\boldsymbol{T},$ at the nominal distribution $F$, is defined as
\begin{align*}
\boldsymbol{\operatorname{IF}}\left(\boldsymbol{z};\boldsymbol{T},F\right)&=\lim_{\epsilon\rightarrow0^+} \frac{\boldsymbol{T}((1-\epsilon)F + \epsilon\Delta_{\boldsymbol{z}})-\boldsymbol{T}(F)}{\epsilon}
\\&=\frac{\partial}{\partial\epsilon} \boldsymbol{T}((1-\epsilon)F + \epsilon\Delta_{\boldsymbol{z}})|_{\epsilon=0},
\end{align*}
where $\epsilon$ is the proportion of samples that are a point-mass, $\Delta_{\boldsymbol{z}},$ located at $\boldsymbol{z}.$

\begin{theorem}[Influence function]
Assuming (A1) and (A2), the influence functions of for the S-q estimates of 
$\widehat{\boldsymbol{\mu}}$ and $\widehat{\boldsymbol{\Sigma}}$
are given by
\begin{align}
    \boldsymbol{\operatorname{IF}}\left(\boldsymbol{z};\boldsymbol{\mu},F\right)
    &=
    \frac{\sqrt{d_z}w_q(d_z/\sigma )}
    {\omega_2}
    \frac{\boldsymbol{z}_c}{\sqrt{d_z}}, \nonumber
    \\ \label{eq:IFSScatter}
    \boldsymbol{\operatorname{IF}}\left( \boldsymbol{z};\boldsymbol{\Sigma},F\right)
    &=\frac{\rho_q\left(d_z/\sigma \right)-b}{\lambda_2}\boldsymbol{\Sigma}
    +\frac{p(p+2)\left(d_z/\sigma\right)w_q\left(d_z/\sigma\right)}{\lambda_1}\left(\frac{\boldsymbol{z}_c\boldsymbol{z}_c^{\T}}{d_z}-\frac{1}{p}\boldsymbol{\Sigma}\right),
\end{align}
where $\boldsymbol{z}_c=\boldsymbol{z}-\boldsymbol{\mu}$ and $d_z= \boldsymbol{z}_c^{\T}\boldsymbol{\Sigma}^{-1}\boldsymbol{z}_c,$ and where the scalars $\omega_2$, $\lambda_1$, and $\lambda_2$ were defined in Theorem \ref{th:asymptoticnormality}.
\end{theorem}
\begin{proof}
See \citep[Corollary 5.2]{Lopuhaa1989OnCovariance} and \citep[Remark 2]{Lopuhaa1997AsymptoticCovariance}.
\end{proof}

\begin{figure} 
\begin{center}
\includegraphics[width=3in]{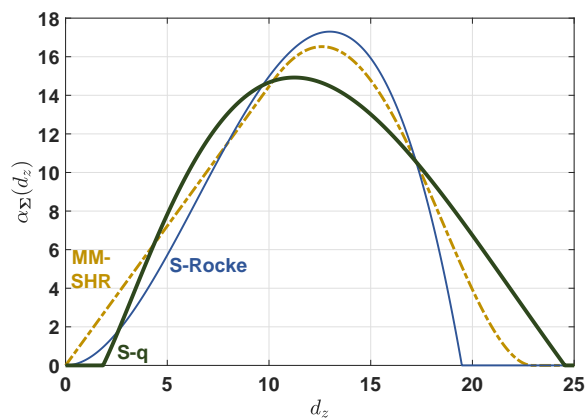}
\end{center}
\caption{Example Comparison of Influence Function Parameter $\alpha_{\boldsymbol{\Sigma}}(d_z).$ For the 10-dimensional Gaussian distribution, $\alpha_{\boldsymbol{\Sigma}}(d_z)$ is plotted for estimators depicted in Figure \ref{fig:compareWs}.
\label{fig:compareIFs}}
\end{figure}

By definition of S-estimators with normalized rho function, the magnitude of first term of (\ref{eq:IFSScatter}) is clearly bounded to no more than $\lambda_2^{-1}\boldsymbol{\Sigma}.$ Therefore, to compare the influence functions of the S-q, S-Rocke, and MM-SHR estimators, we focus on the second term. From this term, define $\alpha_{\boldsymbol{\Sigma}}\left(d_z\right)=\lambda_1^{-1}p(p+2)(d_z/\sigma)w(d_z/\sigma)$ for each estimator. Figure \ref{fig:compareIFs} plots $\alpha_{\boldsymbol{\Sigma}}(d_z)$ at the 10-dimensional Gaussian distribution for the estimators as depicted in Figure \ref{fig:compareWs}.

By definition, all highly-robust estimators have bounded influence functions, and for the three estimators considered here, their influence functions are continuous. This means that small amounts of contamination have limited effects on their estimates. 
The \emph{gross-error sensitivity} of an estimator is the maximum of $\boldsymbol{\operatorname{IF}}\left(\boldsymbol{z}\right),$ and in this example, the S-q demonstrates a lower gross-error sensitivity than the S-Rocke and MM-SHR estimators. By its definition, the MM-SHR has a inlier rejection point of zero, meaning inliers can negatively influence its estimates. However, proper Type II S-q functions have positive inlier rejection points, which provide robustness against inliers.

Relative to the S-Rocke and MM-SHR estimators, the S-q often has larger outlier rejection points. This is the cost of its generally higher efficiency and ability to reject inliers. However, due to its continuity, the influence near this point is still greatly attenuated.

\subsection{Finite-Sample Robustness} \label{subsec:finitesampleRobustness}
\begin{figure} 
\begin{center}
\begin{subfigure}[b]{0.49\textwidth}
\centering
\includegraphics[width=\textwidth]{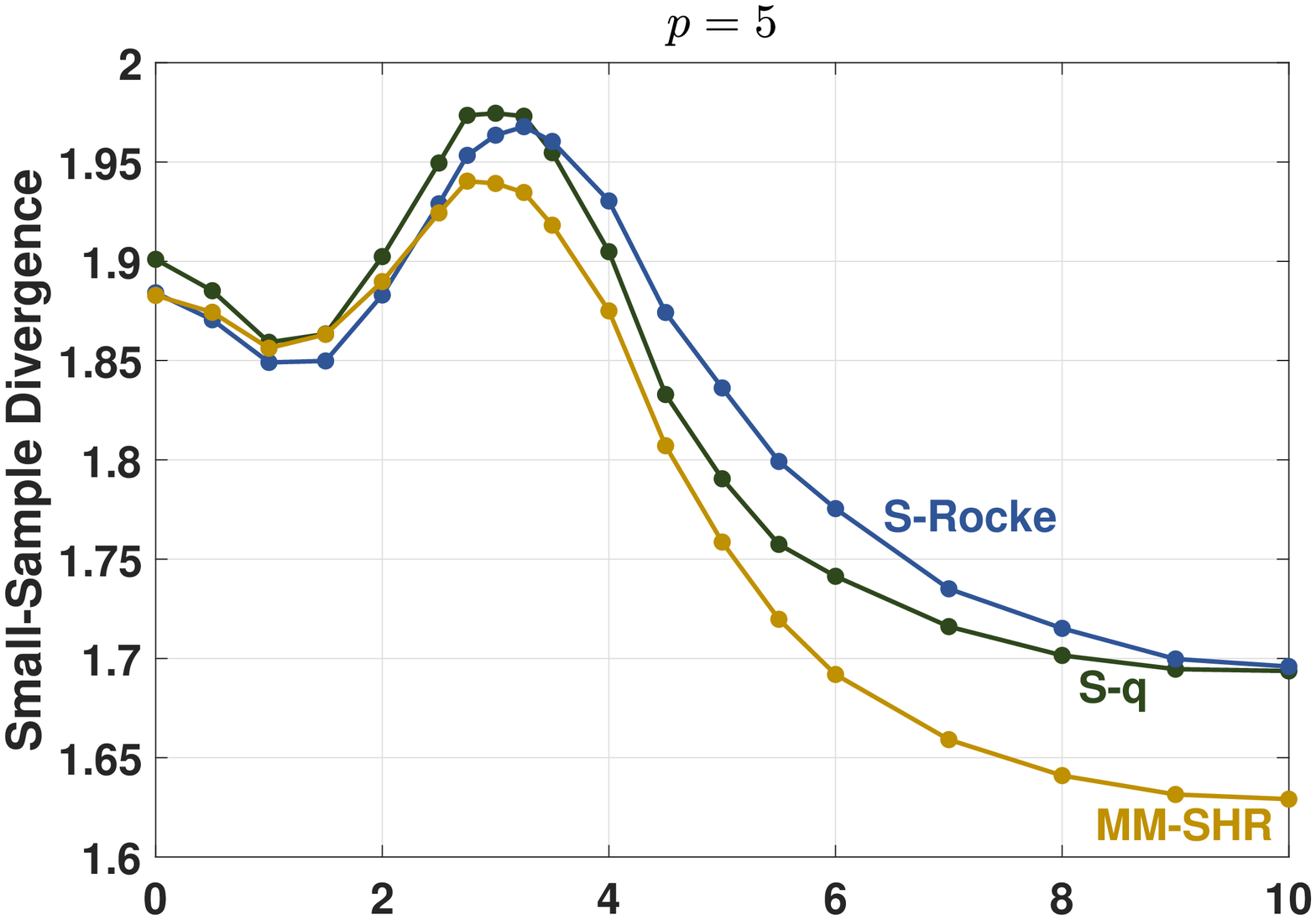}
\end{subfigure}
\hfill
\begin{subfigure}[b]{0.49\textwidth}
\centering
\includegraphics[width=\textwidth]{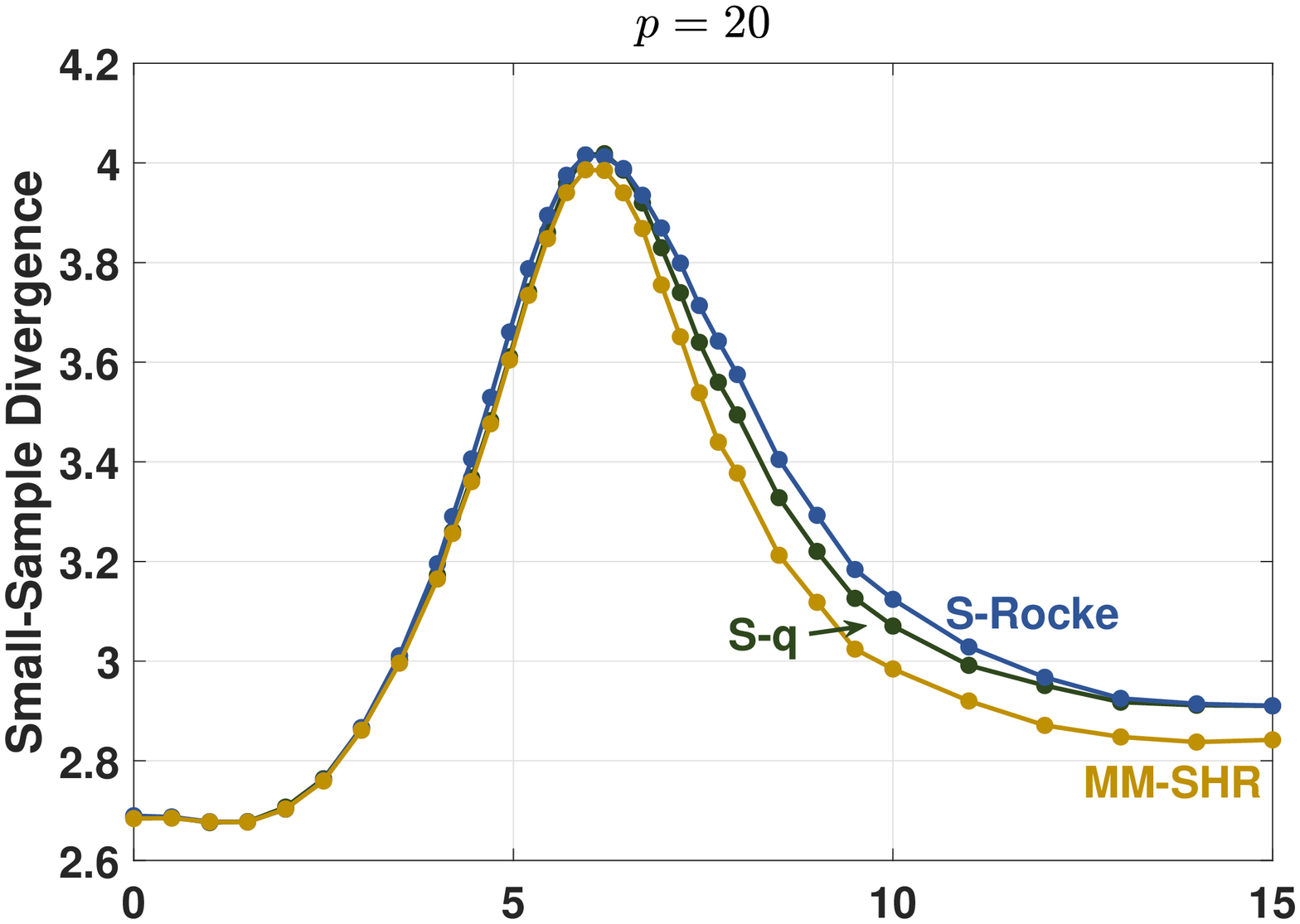}
\end{subfigure}
\newline
\begin{subfigure}[b]{0.49\textwidth}
\centering
\includegraphics[width=\textwidth]{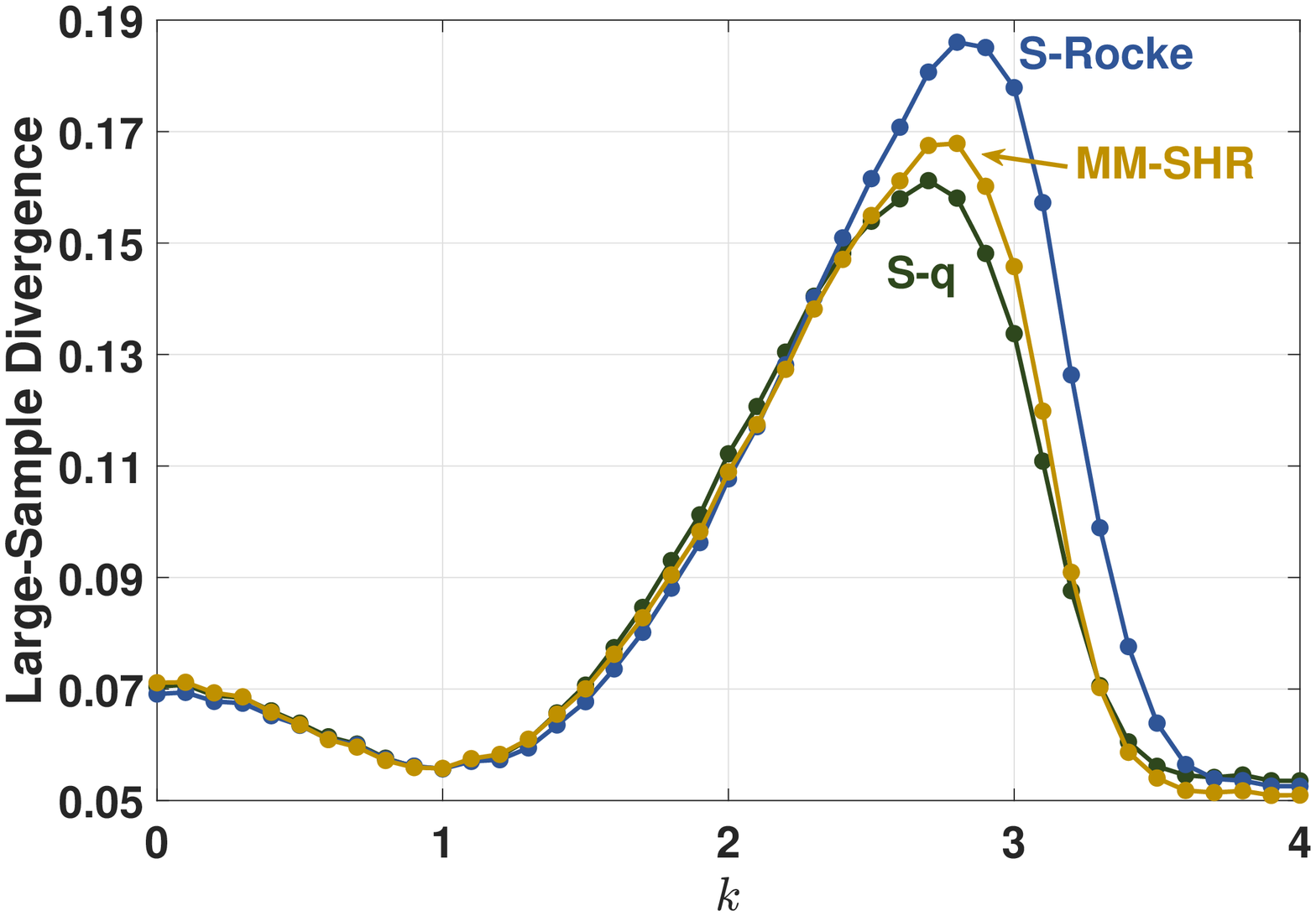}
\end{subfigure}
\hfill
\begin{subfigure}[b]{0.49\textwidth}
\centering
\includegraphics[width=\textwidth]{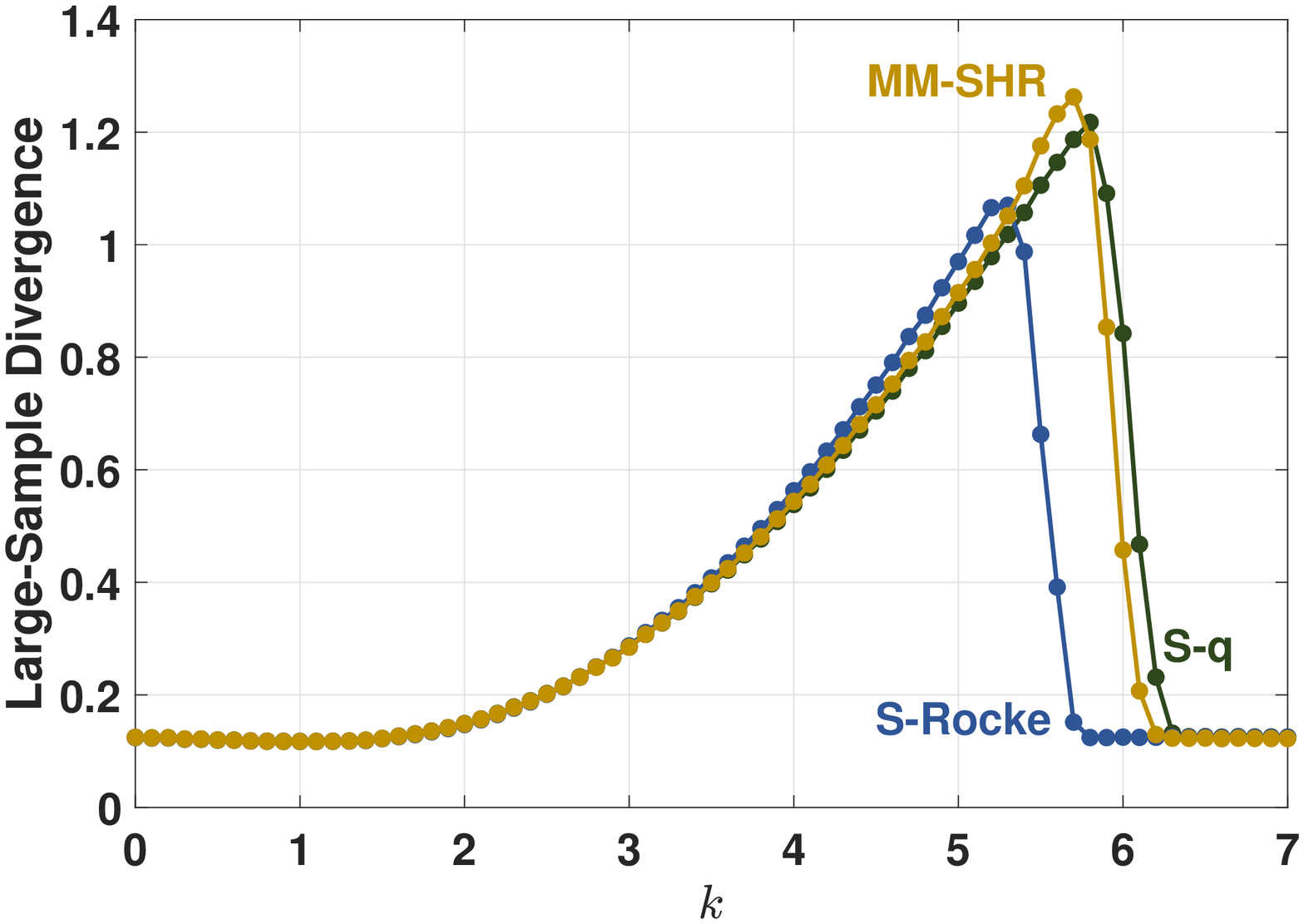}
\end{subfigure}
\end{center}
\caption{Estimator Divergence versus Contamination Value $k$ for Gaussian Distribution. Gaussian shape matrix divergences are plotted for $p=5$ (left) and $p=20$ (right), and for small sample ($n=5p$; top) and large sample ($n=100p$; bottom) sizes.
\label{fig:robustness}}
\end{figure}

To empirically compare the finite-sample robustness of the estimators, we employed the simulation method used by \cite{Maronna2017RobustLocation} and plot the shape matrix divergence, $\D\left(\boldsymbol{\mu},\boldsymbol{\Omega};\boldsymbol{\mu},\widehat{\boldsymbol{\Omega}}\right)$, versus shift contamination value $k$. For a contamination proportion $\epsilon$, the first element of each of the $\lfloor\epsilon n\rfloor$ contaminated samples was replaced with the value $k$, that is $x_1=k$. The initial estimates of the weighted algorithm were determined with the KSD estimator. Figure \ref{fig:robustness} provides divergence plots for normally distributed data with $\epsilon=10\%$ contamination, for dimensions $p=5$ and $p=20$, and for sample sizes $n=5p$ and $n=100p$. For the cases where $p=20$, the estimators were tuned to 90\% uncontaminated relative efficiency. When $p=5$, the S-Rocke has poor maximum efficiency, so the estimators were tuned to match the maximum S-Rocke efficiency.

These plots show that the robustness of the S-q is on par with the other two estimators. Consistent with the results in \cite{Maronna2017RobustLocation}, the relative worst-case performance of the estimators vary by such factors as dimension, sample size, and contamination percentage. For example, the S-q performs the best here for $p=5$, $n=100p$, but the MM-SHR is the best for $p=5$, $n=5p$.

\section{Computational Aspects}
\label{sec:compAnalysis}
This section explores computational aspects of the S-q and other high-breakdown estimators when using their weighted-sum algorithms. The estimators' stabilities are first assessed by comparing their sensitivities to the initial estimates, $\widehat{\boldsymbol{\mu}}^{(0)}$ and $\widehat{\boldsymbol{\Omega}}^{(0)}.$ The computational efficiency is then evaluated by comparing the computational convergence rates of the estimators. 

\subsection{Stability}
The primary criticism of high-breakdown estimators is that their solutions are highly sensitive to the initial estimates $\widehat{\boldsymbol{\mu}}^{(0)}$ and $\widehat{\boldsymbol{\Omega}}^{(0)}$ due to the non-convexity of their objective functions. The S-q estimator helps mitigate this with a generally wider weight function (see for example Figure \ref{fig:compareWs}). To demonstrate that the S-q is more stable with respect to the initial estimates, $m$ Gaussian Monte Carlo simulation trials were run where for each trial, $\boldsymbol{\Omega}$ was estimated twice for each estimator using different initializations. For the first estimate, $\widehat{\boldsymbol{\Omega}}^{(0)}$ was set to the MLE using all $n$ samples before contamination. For the second estimate, $\widehat{\boldsymbol{\Omega}}^{(0)}$ was set to the MLE using just 25\% of the samples before contamination, resulting in a larger expected variance of $\widehat{\boldsymbol{\Omega}}^{(0)}.$ The sample mean of the divergence between the two final estimates, $\bar{\D}\left(\widehat{\boldsymbol{\Omega}}_1,\widehat{\boldsymbol{\Omega}}_2\right),$ was then calculated. The same values of $p,$ $n,$ and $\epsilon$ were used as in the Section \ref{subsec:finitesampleRobustness} simulations. The contamination method was also the same, and the value of $k$ was set to the worst-case value for that estimator and for the values of $p,$ $n,$ and $\epsilon$ (see Figure \ref{fig:robustness}). 

The results are presented in the center of Table \ref{tab:compAnalysis}. It is seen that the S-q estimator was consistently the most stable of the three estimators, and the MM-SHR was generally the most sensitive. For the near-asymptotic cases, the S-q exhibited no measurable differences between the two estimates, unlike the MM-SHR. Like the S-q, the S-Rocke had no measurable differences between the two estimates for the uncontaminated near-asymptotic cases, but under contamination, its mean divergence was roughly on par with the MM-SHR.

\begin{table} \caption{Estimator Stability and Computational Efficiency for Normally Distributed Data} \label{tab:compAnalysis}
\begin{center}
\small
\renewcommand{\arraystretch}{1}
\begin{tabular} {r r r r r r r r r r}
\hline
Dim. & Samples & Contam. & \multicolumn{3}{c}{Mean Divergence} && \multicolumn{3}{c}{Median No. of Iterations}\\
\cline{4-6} \cline{8-10}
$p$ &               $n$ &    $\epsilon$ & S-q & S-Rocke & MM-SHR && S-q & S-Rocke & MM-SHR\\
\hline
        5 &      $100p$ &           0\% &    0 &    0 & 8e-5 &&  13 &     14 &      14\\
        5 &      $100p$ &          10\% &    0 & 7e-4 & 5e-4 &&  14 &     15 &      15\\
       20 &      $100p$ &           0\% &    0 &    0 & 5e-5 &&   7 &      7 &       7\\
       20 &      $100p$ &          10\% &    0 & 1e-3 & 3e-2 &&  18 &      8 &      16\\
        5 &        $5p$ &           0\% & 2e-1 & 4e-1 & 2e-0 &&  32 &     14 &      15\\
        5 &        $5p$ &          10\% & 3e-1 & 5e-1 & 2e-0 &&  31 &     15 &      15\\
       20 &        $5p$ &           0\% & 6e-5 & 4e-2 & 1e-0 &&  28 &     18 &      37\\
       20 &        $5p$ &          10\% & 1e-0 & 2e-0 & 3e-0 &&  30 &     18 &      33\\
\hline
\multicolumn{10}{p{5.74in}}{Note: Mean divergence values listed as ``0'' have simulated average divergences less than the numerical convergence criterion, $\D\left(\widehat{\boldsymbol{\Omega}}^{(i)},\widehat{\boldsymbol{\Omega}}^{(i-1)}\right)<10^{-10}.$}
\end{tabular}
\end{center}
\end{table}

\subsection{Computational Efficiency}
To compare the relative computational efficiencies of the high-breakdown estimators, we calculated the median number of iteration required for the estimators to converge for normally distributed data for various values of $p$, $n$, and $\epsilon$. All three estimators were set to use the same tight convergence criteria that $\D\left(\widehat{\boldsymbol{\Omega}}^{(i)},\widehat{\boldsymbol{\Omega}}^{(i-1)}\right)<10^{-10}.$ The initial estimates were determined with the KSD estimator, and the estimators were tuned as in the Section \ref{subsec:finitesampleRobustness} simulations. The contamination method was also the same, and the value of $k$ was set to the worst-case value for that estimator.

The results are presented on the right of Table \ref{tab:compAnalysis}. For the large-sample ($n=100p$) simulations, the S-q converges approximately as fast as the other two estimators (except for the one case where $p=20,$ $\epsilon=10\%$, where the S-Rocke performs notably better). The S-Rocke estimator consistently converges fastest for all of the small-sample ($n=5p$) cases, and the small-sample convergence of the S-q estimator is relatively consistent---albeit at the upper-end of the spectrum. The small-sample convergence of the MM-SHR is on par with the S-Rocke for small $p,$ but worse than the others for large $p.$

\section{Application to Financial Portfolio Optimization}
\label{sec:example}
A common financial application of mean and covariance matrices is in modern portfolio theory for the optimal allocation of portfolio investments. Under modern portfolio theory's mean-variance framework, a minimum-variance portfolio aims to minimize the risk (i.e. variance) of the portfolio return subject to a desired expected return \citep{Markowitz1952PortfolioSelection}. Mathematically, this is expressed as
\begin{equation*}
\begin{aligned}
\min_{\boldsymbol{\alpha}}
&&& \boldsymbol{\alpha}^{\T}\boldsymbol{\Omega}_r\boldsymbol{\alpha} \\
\text{subject to}
&&& \boldsymbol{\alpha}^{\T}\boldsymbol{\mu}_r=\mu_p,\;\boldsymbol{\alpha}^{\T}\boldsymbol{1}=1,
\end{aligned}
\end{equation*}
where $\boldsymbol{\alpha}$ is a normalized vector of portfolio allocation for each asset, $\boldsymbol{\Omega}_r$ is the shape (or equivalently covariance) matrix for the asset returns, $\boldsymbol{\mu}_r$ is the expected returns of each asset, $\mu_p$ is the desired expected portfolio return, and $\boldsymbol{1}$ is a vector of ones. The solution is given by \citep{Roy1952SafetyAssets,Merton1972AnFrontier}
\begin{align} \label{eq:portfolioAlpha}
\boldsymbol{\alpha} = {} \nonumber
    &   s_r\left(\boldsymbol{\mu}_r^{\T}\boldsymbol{\Omega}_r^{-1}\boldsymbol{\mu}_r\right) \boldsymbol{\Omega}_r^{-1}\boldsymbol{1}
      - s_r\left(\boldsymbol{1}^{\T}\boldsymbol{\Omega}_r^{-1}\boldsymbol{\mu}_r\right) \boldsymbol{\Omega}_r^{-1}\boldsymbol{\mu}_r\\
    & + s_r\mu_p\left(\boldsymbol{1}^{\T}\boldsymbol{\Omega}_r^{-1}\boldsymbol{1}\right) \boldsymbol{\Omega}_r^{-1}\boldsymbol{\mu}_r
      - s_r\mu_p\left(\boldsymbol{\mu}_r^{\T}\boldsymbol{\Omega}_r^{-1}\boldsymbol{1}\right) \boldsymbol{\Omega}_r^{-1}\boldsymbol{1},
\end{align}
where $s_r$ is a scalar that ensures the elements of $\boldsymbol{\alpha}$ sum to one.

In this section, the performances of the MM-SHR, S-Rocke, and S-q estimators are compared for the optimal allocation of investment in the component stocks of the DOW Jones Industrial Average. For each estimator, the parameters $\boldsymbol{\Omega}_r$ and $\boldsymbol{\mu}_r$ were estimated for the daily returns from the component stocks. Then, using a desired portfolio daily return of $\mu_p=.038\%$ (corresponding to 10\% annual return), the optimal allocations, $\boldsymbol{\alpha},$ were calculated using (\ref{eq:portfolioAlpha}). Using $\boldsymbol{\alpha}$ for each estimator, the portfolio return was then calculated for each business day of the verification period, assuming a daily re-balance of investments. Finally, each estimator’s performance was characterized by the variance of these daily returns. This variance is a measure of the volatility of the portfolio.

For the S-q estimator, we noted that \cite{KonlackSocgnia2014AReturns} showed that the generalized hyperbolic distribution is a good model for stock returns, and specifically the variance gamma subclass has good parameter stability over time. Although their analysis is for log returns, daily log returns are generally close to one, so the variance gamma model should also fit well for gross (i.e. linear) returns. For the variance gamma S-q estimator, a density-weighted M-estimator was used to estimate the model parameters $\lambda$ and $\psi$.

\begin{figure} 
\begin{center}
\includegraphics[width=3in]{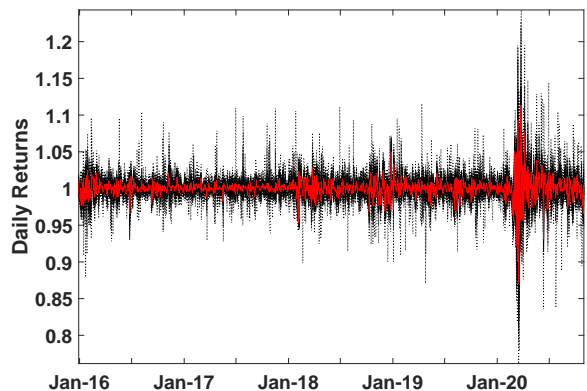}
\end{center}
\caption{Daily Returns of Down Jones Industrial Average and Component Stocks for 2016-2020. There are 25 component stocks in the index throughout the depicted period and plotted in the background. The overall index value is plotted on top.   
\label{fig:stockReturns}}
\end{figure}

To demonstrate the robustness of the S-q estimator, we begin by noting that the first quarter of 2020 contained a once-in-a-generation period of extremely high volatility due to the COVID-19 pandemic, as depicted in Figure \ref{fig:stockReturns}. This volatility started on approximately February 21. Each estimator's performance was assessed by estimating the parameters $\boldsymbol{\Omega}_r$ and $\boldsymbol{\mu}_r$ using all the returns from the first quarter, then comparing the variances of the daily portfolio returns for only the pre-pandemic (prior to February 21) period. Each estimator was set to its maximum breakdown point. Each estimator was then tuned it to its maximum asymptotic efficiency with respect to the variance gamma distribution with parameters estimated using a maximum likelihood approach and using the daily returns for the years 2016--2019. 

Table \ref{tab:varCrash} summarizes the results, listing the variances of the daily returns. The S-q estimator performed the best with the lowest variance, which indicates high robustness. The MM-SHR performed the second best, followed by the S-Rocke. The sample estimator of mean and covariance was also included to demonstrate its poor robustness.

\begin{table} \caption{Sample Variances of Achieved Daily Returns for 01-Jan-2020 -- 20-Feb-2020} \label{tab:varCrash}
\begin{center}
\small
\renewcommand{\arraystretch}{1}
\begin{tabular}{l r}
\hline
Estimator & Variance \\\hline
S-q       & 76 \\
MM-SHR    & 119\\
S-Rocke   & 147\\
Sample    & 176\\
\hline
\end{tabular}
\end{center}
\end{table}

Next, to demonstrate estimator efficiency, variances of daily returns were compared for a non-volatile period: 2016 through 2019. Using the same methodology and configuration as before, for each year and each estimator, $\boldsymbol{\Omega}_r$ and $\boldsymbol{\mu}_r$ were estimated. Then, $\boldsymbol{\alpha}$ was calculated and applied to each day of that year. The sample variances of each year's daily portfolio returns are listed in Table \ref{tab:varYearly}. The S-q estimator resulted in the lowest portfolio variance for three of the four years, and the lowest variance on average, indicating high estimator efficiency. On average, the performance of the MM-SHR estimator was behind that of the S-q estimator, and the S-Rocke demonstrated substantially worse performance.

\begin{table} \caption{Sample Variances of Achieved Daily Returns by Year} \label{tab:varYearly}
\begin{center}
\small
\renewcommand{\arraystretch}{1}
\begin{tabular} {l c c c}
\hline
Year        &          S-q  &       MM-SHR  & S-Rocke \\
\hline
2016        & \textbf{3.51} &         3.56  & 4.10 \\
2017        & \textbf{1.18} &         1.24  & 1.34 \\
2018        & \textbf{6.76} &         6.91  & 7.31 \\
2019        &         3.60  & \textbf{3.50} & 4.49 \\
\hline
Sample Mean & \textbf{3.77} &         3.80  & 4.31 \\
\hline
\end{tabular}
\end{center}
\end{table}

\section{Conclusion}
\label{sec:conclusion}
The S-q estimator has been introduced as a new tunable multivariate estimator of location, scatter, and shape matrices for elliptical probability distributions. This new estimator is a subclass of S-estimators, which achieve the maximum theoretical breakdown point. The S-q estimator has been compared with the leading high-breakdown estimators. Across elliptical distributions, the S-q has generally higher efficiency and stability, and its robustness is on par with these other leading estimators. Additionally, the S-q provides a monotonic and upper-bounded efficiency tuning parameter, which provides simpler tuning than the MM-SHR. The S-q is therefore a broadly applicable estimator, providing practitioners with a good general high-breakdown multivariate estimator that can be used across a broad range of practical applications, such as the optimal portfolio example.

\bibliographystyle{cfg/elsarticle-harv} 
\bibliography{cfg/references}





\end{document}